\documentclass[11pt, reqno]{amsart}
\usepackage{amsmath, amsthm, a4, latexsym, amssymb}

\def\@rmrk#1#2{\refstepcounter
    {#1}\@ifnextchar[{\@yrmrk{#1}{#2}}{\@xrmrk{#1}{#2}}}

\makeatletter\@addtoreset{equation}{section}\makeatother

\newfont{\bfit}{cmbxti10 scaled 2000}
\newfont{\biggi}{cmr12 scaled 2000}
\newtheorem{step}{STEP}

\newcommand{\bes}{\begin{step}}
\newcommand{\es}{\end{step}}
 \newcommand{\DeltaD}{\Delta\!^{\mbox{\tiny d}}}
 \newcommand{\lambdaD}{\lambda^{\mbox{\tiny d}}}

 \newcommand{\eps}{\varepsilon}
 
 \newcommand{\supp}{{\rm supp}\,}

 \newcommand{\esssup}{{\rm esssup}\,}
 \newcommand{\essinf}{{\rm essinf}\,}
 \newcommand{\R}{\mathbb{R}}
 \newcommand{\Z}{\mathbb{Z}}
 \newcommand{\N}{\mathbb{N}}

 \newcommand{\Sym}{\mathfrak{S}}
 \newcommand{\prob}{\mathbb{P}}
 \newcommand{\Prob}{{\rm Prob }}

 \newcommand{\me}{\mathbb{E}}
 
 \renewcommand{\P}{\mathbb{P}}
 \newcommand{\E}{\mathbb{E}}
 \newcommand{\one}{\1}

 \newcommand{\skric}{{\mathcal C}}

 \newcommand{\skrih}{{\mathcal H}}

 \newcommand{\skril}{{\mathcal L}}
 
 \newcommand{\Mcal}{{\mathcal M}}
 \newcommand{\Hcal}{{\mathcal H}}
 \newcommand{\Lcal}{{\mathcal L}}
 \newcommand{\Ocal}{{\mathcal O}}
 \newcommand{\Ccal}{{\mathcal C}}
 \newcommand{\Rcal}{{\mathcal R}}

 \newcommand{\heap}[2]{\genfrac{}{}{0pt}{}{#1}{#2}}
 \newcommand{\sfrac}[2]{\mbox{$\frac{#1}{#2}$}}
 
 \newcommand{\ssup}[1] {{\scriptscriptstyle{({#1}})}}
\def\1{{\mathchoice {1\mskip-4mu\mathrm l}      
{1\mskip-4mu\mathrm l}
{1\mskip-4.5mu\mathrm l} {1\mskip-5mu\mathrm l}}}

\newcommand{\eq}{\begin{equation}}
\newcommand{\en}{\end{equation}}

\newenvironment{Proof}[1]
{\vskip0.1cm\noindent{\bf #1}{\hspace*{0.3cm}}}{\vspace{0.15cm}}

\renewcommand{\subsection}{\secdef \subsct\sbsect}
\newcommand{\subsct}[2][default]{\refstepcounter{subsection}
\vspace{0.15cm}
{\flushleft\bf \arabic{section}.\arabic{subsection}~\bf #1  }
\nopagebreak\nopagebreak}
\newcommand{\sbsect}[1]{\vspace{0.1cm}\noindent
{\bf #1}\vspace{0.1cm}}

{\nopagebreak {\hfill{$\diamond$}}\\ }

\newtheorem{theorem}{Theorem}[section]
\newtheorem{lemma}[theorem]{Lemma}

\newtheorem{prop}[theorem]{Proposition}

\theoremstyle{definition}
\newtheorem{remark}[theorem]{Remark}

\newtheoremstyle{thm}{1.5ex}{1.5ex}{\itshape\rmfamily}{}
{\bfseries\rmfamily}{}{2ex}{}

\newtheoremstyle{rem}{1.3ex}{1.3ex}{\rmfamily}{}
{\itshape\rmfamily}{}{1.5ex}{}
\theoremstyle{rem}
\refstepcounter{subsubsection}

\def\thebibliography#1{\section*{Bibliography}
  \list%
  {\arabic{enumi}.}
    {\settowidth\labelwidth{[#1]}\leftmargin\labelwidth
    \advance\leftmargin\labelsep
    \parsep0pt\itemsep0pt
    \usecounter{enumi}}
    \def\newblock{\hskip .11em plus .33em minus .07em}
    \sloppy                   
    \sfcode`\.=1000\relax}

\def\P{\prob}

\newcommand{\sss}{\scriptscriptstyle}
\newcommand{\sM}{{\sss M}}
\newcommand{\sR}{{\sss R}}

\begin{document}
\title[The universality classes in the parabolic Anderson model]{\Large The universality classes in the parabolic
Anderson model}
\author[Remco van der Hofstad, Wolfgang K\"onig, and Peter M\"{o}rters]{}
\maketitle
\thispagestyle{empty}
\vspace{-0.5cm}

\centerline{\sc By Remco van der Hofstad, Wolfgang K\"onig, and Peter M\"{o}rters}
\renewcommand{\thefootnote}{}
\renewcommand{\thefootnote}{1}

\vspace{-0.5cm}
\centerline{\textit{Eindhoven University of Technology,
Universit\"{a}t Leipzig, and University of Bath}}

\begin{quote}{\small }{\bf Abstract.} We discuss the long time behaviour
of the parabolic Anderson model, the Cauchy problem for the heat equation
with random potential on $\Z^d$. We consider general i.i.d.~potentials and show that
exactly \emph{four} qualitatively different types of intermittent behaviour can occur.
These four universality classes depend on the upper tail of the potential distribution:
(1) tails at $\infty$ that are thicker than the double-exponential tails, (2) double-exponential
tails at $\infty$ studied by G\"artner and Molchanov, (3) a new class called \emph{almost bounded
potentials}, and (4) potentials bounded from above studied by Biskup and K\"onig. The new class (3),
which contains both unbounded and bounded potentials, is studied in both the annealed and the quenched
setting. We show that intermittency occurs on unboundedly increasing islands whose diameter is slowly
varying in time. The characteristic variational formulas describing the optimal profiles of the potential
and of the solution are solved explicitly by parabolas, respectively, Gaussian densities.
\end{quote}


\section{Introduction and main results}

\subsection{The parabolic Anderson model}~\\
\noindent We consider the continuous solution $v\colon[0,\infty)\times \Z^d \to [0,\infty)$
to the Cauchy problem for the heat equation with random coefficients and localised initial datum,
    \begin{eqnarray}
    \frac{\partial}{\partial t} v(t,z) & = & \DeltaD v(t,z) + \xi(z) v(t,z),
    \qquad \mbox{ for } (t,z)\in(0,\infty)\times \Z^d,\label{PAM}\\
    v(0,z) & = & \one_0(z),\qquad \mbox{ for } z\in\Z^d.\label{PAMinitial}
    \end{eqnarray}
Here $\xi=(\xi(z) \colon z\in \Z^d)$ is an i.i.d.~random potential
with values in $[-\infty,\infty)$, and $\DeltaD$ is the discrete Laplacian,
    $$
    \DeltaD f(z)= \sum_{y \sim z} \bigl[f(y)-f(z)\bigr], \qquad
    \mbox{ for } z\in \Z^d,\, f\colon\Z^d\to\R.
    $$
The parabolic problem~\eqref{PAM} is called the \emph{parabolic Anderson model}. The operator $\DeltaD+\xi$
appearing on the right is called the \emph{Anderson Hamiltonian}; its spectral properties are well-studied in
mathematical physics. Equation~\eqref{PAM} describes a random mass transport through a random field of
sinks and sources, corresponding to lattice points $z$ with $\xi(z)<0$, respectively, $>0$.
It is a linearised model for chemical kinetics \cite{GM90}, is equivalent to
Burger's equation in hydrodynamics \cite{CM94}, and describes magnetic phenomena
\cite{MR94}. We refer the reader to \cite{GM90}, \cite{M94} and \cite{CM94} for more background
and to \cite{GK04} for a survey on mathematical results.

The long-time behaviour of the parabolic Anderson problem is well-studied in the mathematics and
mathematical physics literature because it is the prime example of a model exhibiting an
\emph{intermittency effect}. This means, loosely speaking, that most of the total mass of the solution,
\begin{equation}\label{Utdef}
U(t)=\sum_{z\in\Z^d} v(t,z), \qquad\mbox{for } t>0,
\end{equation}
is concentrated on a small number of remote islands, called the {\it intermittent islands}.
A manifestation of intermittency in terms of the moments of $U(t)$ is as follows. For
$0<p<q$, the main contribution to the $q^{\rm th}$ moment of $U(t)$  comes from islands
that contribute only negligibly to the $p^{\rm th}$ moments. Therefore, intermittency
can be defined by the requirement,
\begin{equation}\label{Intermitt}
\limsup_{t\to\infty}\frac {\langle U(t)^p\rangle^{1/p}}{\langle U(t)^q\rangle^{1/q}}=0,
\qquad \mbox{ for $0<p<q$, }
\end{equation}
where $\langle\,\cdot\,\rangle$ denotes expectation with respect to $\xi$. Whenever $\xi$ is truly
random, the parabolic Anderson model is intermittent in this sense, see \cite[Theorem~3.2]{GM90}.

However, one wishes to understand the intermittent behaviour in much greater detail. The following has been
heuristically argued in the literature and has been verified, at least partially, for important special
examples of potentials: the intermittent islands are characterized by a particularly high exceedance of the
potential and an optimal shape, which is determined by a deterministic variational formula. A universal
picture is present: the location and number of the intermittent islands are random, their size and the absolute
height of the potential in the islands is $t$-dependent, but the (rescaled) shape depends neither on
randomness nor on $t$. Examples studied include the double-exponential distribution
\cite{GM98}, potentials bounded from above \cite{BK01} and continuous analogues
on $\R^d$ instead of $\Z^d$ like Poisson obstacle fields \cite{S98} and Gaussian and other Poisson
fields \cite{GK98, GKM99}. A finer analysis of the geometry of the intermittent islands has been
carried out for Poisson obstacle fields \cite{S98} and the double-exponential
distribution \cite{GKM05}.

In the present paper we initiate the study of the parabolic Anderson model for \emph{arbitrary}
potentials, with the aim of identifying {\it all universality classes\/} of intermittent behaviour
that can arise for different potential distributions. Our standing assumption is that the potentials
$(\xi(z) \colon z\in\Z^d)$ are independent and identically distributed and that all positive exponential
moments of $\xi(0)$ are finite, which is necessary and sufficient
for the finiteness of the $p^{\rm th}$ moments of $U(t)$ at all times. The long-term behaviour of the
solutions depends strongly and exclusively on the upper tail behaviour of the random variable $\xi(0)$.
It is fully described by the top of the spectrum of the Anderson Hamiltonian
$\DeltaD+\xi$ in large $t$-dependent boxes.

The outline of the remainder of this section is as follows. In Section~\ref{sec-reg}, we formulate and
discuss a mild regularity condition on the potential. In Section~\ref{sec-uni}, we show that under this
condition the potentials can be split into exactly \emph{four} classes, which exhibit four different
types of intermittent behaviour.
Three of these classes have been studied in the literature up to now. A fourth class,
the class of \emph{almost bounded potentials}, is studied in the present paper for the first time. We
present our results on the moment and almost-sure large-time asymptotics for $U(t)$ in
Section~\ref{sec-abp}. In Section~\ref{MA}, we give a heuristic derivation of the moment asymptotics,
and in Section~\ref{Varprob}, we explain the variational problems involved.

\subsection{Regularity assumptions}\label{sec-reg}

\noindent We first state and discuss our regularity assumptions on the potential.
Roughly speaking, the purpose of these assumptions is to ensure that
the potential has the same qualitative behaviour at different scales, and therefore
the system does not belong to different universality classes at different times.
Our assumptions refer to the upper tail of $\xi(0)$, and are conveniently formulated
in terms of the regularity of its logarithmic moment generating function,
    \begin{equation}\label{Hdef}
    H(t)= \log \big\langle e^{t\xi(0)} \big\rangle, \qquad\mbox{ as } t \uparrow \infty.
    \end{equation}
Note that $H$ is convex and $t\mapsto H(t)/t$ is increasing with
$\lim_{t\to \infty} H(t)/t=\esssup \xi(0)$.
To simplify the presentation, we make the assumption that if $\xi$ is bounded from above,
then $\esssup \xi(0)=0$, so that $\lim_{t\to \infty} H(t)/t\in \{0,\infty\}$.
This is no loss of generality, as additive constants in the potential appear as additive
constants both in $\frac 1{pt} \log \langle U(t)^p \rangle$ and $\frac 1t \log U(t)$.
The first central assumption on $H$ is the following:

\medskip
{\bf Assumption (H).}
\textit{$t \mapsto \frac{H(t)}{t}$ is in the de Haan class.}
\medskip

We recall that a measurable function $\widetilde{H}$ is said to
be in the \emph{de Haan class} if, for some
regularly varying function $g\colon(0,\infty)\to\R$, the term
$g(t)^{-1}\big(\widetilde{H}(\lambda t)-\widetilde{H}(t)\big)$
converges to a nonzero limit as $t\uparrow \infty$, for any $\lambda>1$.
Recall that a measurable function $g$ is called \emph{regularly varying} if
$g(\lambda t)/g(t)$ converges to a positive limit for every
$\lambda>0$. If this is the case, then the limit takes the form $\lambda^\varrho$,
and $\varrho$ is called the \emph{index of regular variation}. If $\varrho=0$,
then the function is called \emph{slowly varying}.

When $H(t)/t$ is in the de Haan class, then $H$ is regularly
varying with some index $\gamma\in\R$. By convexity of $H$, we have
$\gamma\geq 0$.  If $H$ is regularly varying with index $\gamma\not=1$,
then $H(t)/t$ is in the de Haan class, so that the statements are
equivalent for $\gamma\not=1$. However, if $\gamma=1$, then this
does not necessarily hold,
see \cite[Theorem 3.7.4]{BGT87}.

From the theory of regular functions we derive the existence of a function
$\widehat H$ which can be characterized by two parameters,
$\gamma\in[0,\infty)$ and $\rho\in(0,\infty)$,
and plays an important role in the sequel.

\begin{prop}\label{variation}
Assumption (H) is equivalent to the existence of a function $\widehat{H}\colon(0,\infty)\to\R$ and
a continuous auxiliary function $\kappa\colon(0,\infty)\to(0,\infty)$ such that
    \begin{equation}\label{basic}
    \lim_{t\uparrow\infty} \frac{H(ty)-yH(t)}{\kappa(t)}
    = \widehat{H}(y)\not=0, \qquad\textit{ for } y\in(0,1) \cup (1,\infty).
    \end{equation}
The convergence holds uniformly on every interval $[0,M]$, with~$M>0$.
Moreover, with $\gamma$ the index of variation of $H$,
the following statements hold:
\begin{itemize}
\item[(i)] $\kappa$ is regularly varying of index $\gamma\geq 0$. In particular,
$\kappa(t)=t^{\gamma+o(1)}$ as $t\uparrow\infty$.
\item[(ii)] There exists a parameter $\rho >0$ such that, for every $y>0$,\\[-1mm]
\begin{itemize}
\item[(a)] if $\gamma\not= 1$, then $
\displaystyle\widehat H(y)=\rho \, \frac{y-y^\gamma}{1-\gamma}$,
and $\displaystyle\lim_{t\uparrow\infty}\frac{H(t)}
{\kappa(t)}=\frac\rho{\gamma-1},$ \\[1mm]
\item[(b)] if $\gamma=1$, then $\widehat H(y)=\rho  y \log y$, and
$\displaystyle\lim_{t\uparrow\infty}\frac{| H(t) | }
{\kappa(t)}=\infty$.
\end{itemize}
\end{itemize}
\end{prop}

\begin{Proof}{Proof.}
See Chapter~3 in \cite{BGT87}. More accurately, using the notation
$f(t)=H(t)/t$ and $g(t)=\kappa(t)/t$, (i) is shown in \cite[Section~3.0]{BGT87},
see also \cite[Theorem~1.4.1]{BGT87}. The uniformity of the convergence follows since
the left hand side of \eqref{basic}
is convex in $y$, negative on the interval $(0,1)$, and
continuous in zero.

(ii) follows from \cite[Lemma~3.2.1]{BGT87}.
The implication stated in (ii)(a) follows from \cite[Theorems~3.2.6, 3.2.7]{BGT87},
and the implication stated in (ii)(b) is shown in \cite[Theorem~3.7.4]{BGT87}.
\qed
\end{Proof}

Note that $\kappa$ is an asymptotic scale function, and $\widehat H$ an asymptotic
shape function for $H$. While $\gamma\in[0,\infty)$ is unambiguously determined by
the potential distribution, the parameter $\rho$ could be absorbed in either
$\kappa$ or $\widehat H$. The latter option makes it possible to keep track of $\rho$
in the sequel. If $\xi$ is unbounded from above, then $\xi$ and $\xi+ C$ have the same pair
of $\widehat H$ and $\kappa$ for any $C\in\R$. If $\xi$ is replaced by $C\xi$ for some $C>0$,
then the pair $(\widehat H,\kappa) $ may be  replaced by $(C^\gamma \widehat H,\kappa)$.
In the case $\gamma\not=1$ one may choose $\kappa(t)=H(t)$ in \eqref{basic}, if
$\gamma=1$ one may take $\kappa(t) = H(t) - \int_1^t H(s)/s \,ds$, see \cite[Theorem 3.7.3]{BGT87}.

The three regimes $0\leq \gamma<1$, $\gamma=1$ and $\gamma>1$ obviously distinguish three
qualitatively different classes of (upper tail behaviour of) potentials.
However, in order to appropriately describe the asymptotics of the parabolic Anderson model
in the case $\gamma=1$, a finer distinction is necessary. For this we need an additional mild
assumption on the auxiliary function $\kappa$:
\medskip

{\bf Assumption (K).}
\textit{The limit $\displaystyle \kappa^*=\lim_{t \to \infty} \frac{\kappa(t)}{t}$
exists as an element of $[0,\infty]$.}
\medskip

Assumption (K) is obviously satisfied in the cases $\gamma\not=1$ and
for potentials bounded from above in the case $\gamma=1$. Indeed,
when $\gamma<1$, then $\kappa^*=0$, while when $\gamma>1$, then
$\kappa^*=\infty$ by Proposition \ref{variation}(ii)(a). When $\gamma=1$
and $H(t)/t\to 0$, then, by Proposition \ref{variation}(ii)(b),
$H(t)/\kappa(t) \to \infty$, so that $\kappa(t)/t\to 0$.
Hence, Assumption (K) can be a restriction only for potentials unbounded
from above in the case~$\gamma=1$.

\subsection{The universality classes}\label{sec-uni}

\noindent In this section, we define and discuss the four universality classes of the parabolic
Anderson model under the Assumptions (H) and (K). In particular, we explain the relation
between the asymptotics of the parabolic Anderson model and the parameters
$\gamma$ and $\kappa^*$ introduced in Assumptions~(H) and~(K).

For the moment, we focus on the large time behaviour of the $p^{\rm th}$ moment $\langle U(t)^p \rangle$
for any $p>0$. We show that there is a scale function $\alpha\colon(0,\infty)\to(0,\infty)$ and a number 
$\chi\in \R$ such that
\begin{equation}\label{generalform}
\frac 1{pt}\log\langle U(t)^p \rangle = \frac{H\big(pt \,\alpha(pt)^{-d}\big)}{pt\,\alpha(pt)^{-d}} -
 \frac{1}{\alpha(pt)^2} \, \big( \chi + o(1) \big),
\qquad \mbox{ as } t\uparrow\infty,
\end{equation}
The scale function $\alpha$ describes how fast
the expected total mass, which at time $t=0$ is localised at the origin, spreads, in the sense that
\begin{equation}\label{UFK2}
\begin{aligned}
\lim_{R\uparrow\infty} \liminf_{t\uparrow\infty}
\frac{\alpha(t)^2}{t} \log \, \frac{\left\langle \sum_{z\in \Z^d}
v(t,z)\, \one\{ |z|\leq R\,\alpha(t) \} \right\rangle}
{\left\langle \sum_{z\in \Z^d} v(t,z)\right\rangle} = 0 .
\end{aligned}
\end{equation}
\emph{Heuristically}, $\alpha(t)$ also determines the size of the
intermittent islands for the almost sure behaviour of $U(t)$. The order
of their diameter is given as $(\alpha \circ \beta)(t)$, where $\beta(t)$
is the asymptotic inverse of $t\mapsto t/\alpha(t)^2$ evaluated at $d\log t$, 
cf.~Section~\ref{sec-asasym} below. The numbers $\chi$ are naturally given in 
terms of minimisation problems, where the minimisers correspond to the typical shape
of the solution on an intermittent island. A rigorous proof of these heuristic
statements, however, is beyond the means of this paper.

One expects that $\alpha(t)$ is asymptotically the larger, the thinner the upper tails
of $\xi(0)$ are. It will turn out that when $\kappa^*=\infty$, then
(\ref{UFK2}) is satisfied with $\alpha(t)=1$ independently of $R$.
Therefore, we only need to analyse $\alpha(t)$ in the case when
$\kappa^*<\infty$. Analytically, if $\kappa^*<\infty$, then $\alpha(t)$
may be defined by a fixed point equation as follows:

\begin{prop}[The scale function $\alpha$]
\label{pro-exisalpha}
Suppose that Assumptions~(H) and (K) are satisfied and $\kappa^*<\infty$. There exists a
regularly varying scale function $\alpha\colon(0,\infty)\to(0,\infty)$, which is unique up
to asymptotic equivalence, such that for all sufficiently large $t>0$
 \begin{equation}\label{alphadef}
 \frac{\kappa\bigl(t\alpha(t)^{-d}\bigr)}{t\alpha(t)^{-d}} = \frac{1}{\alpha(t)^{2}}.
 \end{equation}
The index of regular variation is $\frac{1-\gamma}{d+2-d\gamma}$ and hence
$\displaystyle \lim_{t\uparrow\infty} \frac{t}{\alpha(t)^d}=\infty.$ Moreover,
\begin{itemize}
\item[(i)] If $\gamma=1$ and $0<\kappa^*<\infty$, then $\lim_{t\uparrow\infty} \alpha(t)=
1/\sqrt{\kappa^*}\in(0,\infty)$.
\item[(ii)] If $\gamma=1$ and $\kappa^*=0$, or if  $\gamma<1$, then
$\lim_{t\uparrow\infty} \alpha(t)=\infty$.
\end{itemize}
\end{prop}


\begin{Proof}{Proof.} To see that $\alpha$ is regularly varying and unique up to asymptotic
equivalence we note that $f(t)=t (\kappa(t)/t)^{-d/2}$ is regularly varying with index at least one. By
\cite[Theorem 1.5.12]{BGT87}, there exists an asymptotically unique inverse $g$ such
that $f(g(t)) \sim t$ for $t\uparrow \infty$. This inverse is regularly varying.
By definition, $t \mapsto t \alpha(t)^{-d}$ satisfies $f(t \alpha(t)^{-d})=t$ and
hence $\alpha(t) \sim (t/g(t))^{1/d}$ is regularly varying. The  index of regular variation
of $\alpha$ is immediate from the defining equation and the fact that $\kappa(t)$ is regularly 
varying with index $\gamma$.

Under the assumptions of~(i), for large $t$, the mapping $x \mapsto \kappa(tx^{d/2})/tx^{d/2}$ maps
a compact interval centred in $\kappa^*$ to itself, and hence the existence of a solution to~\eqref{alphadef}
follows from a fixed-point argument. The stated properties of $\alpha(\,\cdot\,)$ follow immediately
from the definition.

Under the assumptions of~(ii), we look at the problem of finding $s>0$ such that
$\kappa(s)/s= ({s}/{t})^{2/d}.$
For any fixed $t$, as we increase $s$  the left hand side goes to zero
and the right hand side to infinity. Hence for sufficiently large $t$,
there exists a solution $s=s(t)$, which is going to infinity as $t\uparrow\infty$.
Then $\alpha(t)=(t/s(t))^{1/d}$ solves~\eqref{alphadef} and converges to infinity.
\qed\end{Proof}

Now we introduce the four universality classes, ordered from thick to thin upper tails of $\xi(0)$.
Recall the general formula for the 
asymptotics of the moments $\langle U(t)^p\rangle$ from~\eqref{generalform}.

\begin{itemize}
\item[(1)] \framebox{$\gamma>1$, or $\gamma=1$ and $\kappa^*=\infty$.}\\

\noindent
This case is included in \cite{GM98} as the upper boundary case  $\rho=\infty$ in their notation.
Here $\chi=2d$, the scale function $\alpha(t)=1$ is constant, and the first term on the right hand side
in~\eqref{generalform} dominates the sum, which diverges to infinity.
The asymptotics in \eqref{UFK2} can be strengthened
to $$\lim_{t\uparrow \infty}
\frac 1t \log \frac{\langle v(t,0) \rangle}{ \big\langle \sum_{z \in\Z^d} v(t,z) \big\rangle } = 0,$$
i.e. the expected total mass remains essentially in the origin and the intermittent islands are single
sites, a phenomenon of \emph{complete localisation}.  We call this the {\it single-peak case}.
\hfill$\Diamond$\\

\item[(2)] \framebox{$\gamma=1$ and $\kappa^*\in(0,\infty)$.}\\

\noindent
This case, the \emph{double-exponential case}, is the main objective of \cite{GM98}.
The prime example is the double exponential distribution with parameter $\rho\in(0,\infty)$,
$$
{\rm Prob}\big\{ \xi(0) > r \big\} =\exp\{ - e^{r/\rho}\},
$$
which implies $H(t)= \rho t \log (\rho t) - \rho t + o(t)$.
Here $\alpha(t)\to 1/\sqrt\kappa^*\in(0,\infty)$, so that the size of the
intermittent islands is constant in time. The first term on the right
hand side in \eqref{generalform} dominates the sum, which goes to infinity. Moreover,
\begin{equation}\label{chicase2}
\chi= \min_{\heap{g\colon\Z^d\to\R}{\sum g^2=1}} \Big\{ \frac 12\sum_{\heap{x,y\in\Z^d}{x\sim y}}
\big( g(x)-g(y) \big)^2 - \rho \sum_{x\in\Z^d} g^2(x) \log g^2(x) \Big\},
\end{equation}
where we write $x\sim y$ if $x$ and $y$ are neighbours.
This variational problem is difficult to analyse. It has
a solution, which is unique for sufficiently large values of $\rho$, and
heuristically this minimizer represents the shape of the solution.
As noted in~\cite{GH99}, for any family of minimizers $g_\rho$, as $\rho\uparrow\infty$,
$g_\rho$ converges to $\delta_0$, which links to the single-peak case. Furthermore,
as $\rho\downarrow 0$, the minimisers $g_\rho$ are asymptotically given by
$$
g_\rho^2(\lfloor x/\sqrt\rho\rfloor)=(1+o(1))\, e^{-|x|^2}\pi^{-d/2},
$$ uniformly on
compacts and in $L^1(\R^d)$. Consequently,
$$
\chi=\rho \, d \, \Big(1-\frac 12 \log
\frac\rho\pi+o(1)\, \Big)\qquad \mbox{as } \rho\downarrow 0.\vspace{-0.7cm}
$$\hfill$\Diamond$
\ \\

\pagebreak[2]

\item[(3)] \framebox{$\gamma=1$ and $\kappa^*=0$.}\\

\noindent
Potentials in this class are called \emph{almost bounded} in \cite{GM98} and may be seen as the
degenerate case for $\rho=0$ in their notation. This
class contains both bounded and unbounded potentials, and is analysed for the first time in the
present paper. The scale function $\alpha(t)$ and hence the diameter of the intermittent islands
goes to infinity and is slowly varying, in particular it is slower than any power of $t$. The first term
on the right hand side in \eqref{generalform} dominates the sum, which may go to infinity
or zero. Moreover,
\begin{equation}\label{chicase3}
\chi= \min_{\heap{g\in H^1(\R^d)}{\|g\|_2=1}} \Big\{ \int_{\R^d}|\nabla g(x)|^2\, dx -
\rho \int g^2(x) \log g^2(x) \, dx \Big\},
\end{equation}
see Theorem~\ref{momasy}. This variational formula is obviously the continuous variant
of \eqref{chicase2}, and it is much easier to solve. There is a unique minimiser,
given by
$$
g_\rho(x)=\Big(\frac\rho\pi\Big)^{d/4}\exp\Big(-\frac\rho2 |x|^2\Big),
$$
representing the rescaled shape of the solution on an intermittent island.
In particular, $\chi=\rho d\big(1-\frac 12 \log \sfrac \rho \pi\big)$, which
is the asymptotics of \eqref{chicase2} as $\rho\downarrow 0$. Hence, on the level of
variational problems, (3) is the boundary case of (2) for $\rho\downarrow 0$.\hfill$\Diamond$
\vspace{0.5cm}

\item[(4)] \framebox{$\gamma<1$.}\\

\noindent
This is the case of \emph{potentials bounded from above}, which is treated in \cite{BK01}. Indeed,
in \cite{BK01}, it is assumed that there exists a non-decreasing function $\alpha(t)$
and a nonpositive function $\widetilde{H}\colon (0,\infty) \to (-\infty,0]$ such that
$$\lim_{t \uparrow \infty} \sfrac{\alpha(t)^{d+2}}{t} H\big( \sfrac{t}{\alpha(t)^d}\, y\big)
= \widetilde{H}(y),$$
uniformly on compact sets in $(0,\infty)$. It is easy to infer from the results of Section~\ref{sec-reg}
above that this assumption holds if Assumption (H) holds for the index $\gamma<1$, for $\alpha$ defined
by~\eqref{alphadef} and
    $$
    \widetilde H(y)=\frac\rho{\gamma-1}\,y^\gamma.
    $$ 
Here $\alpha(t)\to\infty$ as $t\mapsto\alpha(t)$ is regularly varying with index
$\frac{1-\gamma}{d+2-d\gamma}$. The potential $\xi$ is necessarily bounded from above.
In this case, the two terms on the right hand side in \eqref{generalform} are of the \emph{same order},
and  \eqref{generalform} converges to zero. Moreover,
    \begin{equation}\label{Bikoe}
    \chi= \inf_{\heap{g\in H^1(\R^d)}{\|g\|_2=1}} \Big\{ \int_{\R^d}|\nabla g(x)|^2\, dx -
    \rho\int_{\R^d}\frac {g^{2\gamma}(x)-g^2(x)}{\gamma-1}\, dx \Big\}.
    \end{equation}
In the lower boundary case where $\gamma=0$, the functional $\int g^{2\gamma}$ must be replaced by
the Lebesgue measure of $\supp(g)$. In this case the formula is well-known and well-understood.
In particular, the minimizer exists, is unique up to spatial shifts,  and has compact support.
To the best of our knowledge, for $\gamma\in(0,1)$, the formula in \eqref{Bikoe} has not
been analysed explicitly, unless in $d=1$.
In Proposition \ref{prop-convgammaVP} below, we show that \eqref{Bikoe} converges to
\eqref{chicase3}, as would follow from interchanging the limit $\gamma\uparrow 1$
with the infimum on $g$. This means that, on the level of variational formulas, (3)
is the boundary case of (4) for $\gamma\uparrow 1$.
\end{itemize}

\begin{remark}\label{rem-chiinter}
The variational problems in \eqref{chicase2}, \eqref{chicase3}, and \eqref{Bikoe}
encode the asymptotic shape of the rescaled and normalised solution $v(t,\,\cdot\,)$ in the centred ball
with radius of order~$\alpha(t)$. Informally, the main contribution to $\langle U(t) \rangle$
comes from the events that $$\frac{v\big(t, \lceil \, \cdot \, \alpha(t) \rceil \big)}
{\big\|v\big(t, \lceil \, \cdot \, \alpha(t) \rceil \big)\big\|_2} \approx g,$$
where $g$ is a minimiser in the definition of $\chi$. To the best of our knowledge this heuristics 
has not been made rigorous in any nontrivial case so far. Note that in case~(1), formally, \eqref{chicase2} 
holds with $\rho=\infty$ and hence the optimal $g$ is $\one_0$.
\hfill$\Diamond$
\end{remark}

Since the cases (1), (2) and and (4) have been studied in the literature
\cite{BK01, GM98}, the possible scaling picture of the parabolic Anderson
model under the Assumptions (H) and (K) is complete once the case (3)
is resolved. This is the content of the remainder of this paper.

\subsection{Long time tails in the almost bounded case}\label{sec-abp}

\noindent In this section we present our results on the almost bounded case~(3).
In other words, we assume that $\kappa(t)/t$ is slowly varying and converges to zero.

\subsubsection{Moment asymptotics}

\noindent Our main result on the annealed asymptotics of $U(t)$ gives the first two terms
in the asymptotics of $\langle U(t)^p \rangle$ for any $p>0$, as $t\uparrow\infty$.

\begin{theorem}[Moment asymptotics]\label{momasy}
Suppose Assumptions (H) and (K) hold, and assume that we are in case~{\rm (3)}, i.e.,
$\gamma=1$ and $\kappa^*=0$. Let $\rho >0$ be as in
Proposition~\ref{variation}(ii)(b). Then, for any $p\in(0,\infty)$,
\begin{equation}\label{result1}
\frac 1{pt}\log\langle U(t)^p \rangle = \frac{H\big(pt \,\alpha(pt)^{-d}\big)}{pt\,\alpha(pt)^{-d}} -
 \frac{1}{\alpha(pt)^2} \, \big(\rho d(1-\sfrac 12 \log \sfrac \rho \pi) +o(1)\big),
\qquad \mbox{ as } t\uparrow\infty.
\end{equation}
\end{theorem}

\begin{remark}[The constant] Recall from \eqref{generalform} and \eqref{chicase3}
that the constant $\rho d(1-\frac 12 \log \sfrac \rho \pi)$ arises as a variational
problem; see Section \ref{Varprob}. The variational problem plays an essential role
in the proof.\hfill$\Diamond$
\end{remark}

\begin{remark}[Intermittency]\label{rem-inter}
Note from \eqref{alphadef} that the first term in \eqref{result1} is of
higher order than the second term. Formula \eqref{result1}, together with the
results of Proposition~\ref{variation} and the fact that $\alpha(\,\cdot\,)$ is slowly varying, imply that
\begin{equation}\label{intermittency}
\begin{aligned}
\log\frac{ \langle U(t)^p\rangle^{1/p}}{
\langle U(t)^q\rangle^{1/q}} =
\frac{H\big(pt \,\alpha(pt)^{-d}\big)}{p\,\alpha(pt)^{-d}} -
\frac{H\big(qt \,\alpha(qt)^{-d}\big)}{q\,\alpha(qt)^{-d}} + o\big(t/\alpha(t)^2\big) \\
= \frac{t}{\alpha(t)^2} \, \Big( \sfrac qp \hat{H}\big(\sfrac pq \big) + o(1) \Big)
\qquad\mbox{ for } p,q\in(0,\infty).
\end{aligned}
\end{equation}
In particular, we have intermittency in the sense of  \eqref{Intermitt}, and the convergence is
exponential on the scale~$t/\alpha(t)^2$.
\hfill$\Diamond$
\end{remark}

In spite of the simplicity of the variational formula \eqref{chicase3}, the derivation
of \eqref{result1} is technically rather involved and requires a number of demanding tools.
We use both representations of $U(t)$ available to us: an approximative representation in terms of an
eigenfunction expansion, and the Feynman-Kac formula involving simple random walk. The heart
of the proof is an application of a large deviation principle for the rescaled
local times of simple random walk. However, there are three major obstacles to be removed,
which require a variety of novel techniques. The first one is a compactification argument
for the space, which is based on an estimate for Dirichlet eigenvalues in large boxes
against maximal Dirichlet eigenvalues in small subboxes. This is an adaptation of a
method from~\cite{BK01}. The second technique is a cutting argument for the large potential
values, which we trace back to a large deviations estimate for the self-intersection number
of the simple random walk. This is of independent interest and is carried out in Section~\ref{sec-cutting}.
Finally, the third obstacle, which appears in the proof of the \emph{upper} bound,
is the lack of upper semi-continuity of the map $f\mapsto \int f(x)\log f(x)\, dx$
in the topology of the large deviation principle, even after compactification and removal of large values.
Therefore, in the proof of the upper bound we replace the classical large deviation principle
by a new approach, taken from \cite{BHK04}, which identifies and estimates the joint density
of the family of the random walk local times. See Proposition \ref{prop-BHK04} below.

An alternative heuristic derivation of formula \eqref{result2} is given in Section~\ref{MA}.
The proof of Theorem~\ref{momasy} is given in Sections~\ref{sec-cutting} and~\ref{sec-proofann}.

\subsubsection{Almost-sure asymptotics}\label{sec-asasym}

We define another scale function $\beta$ such that
\begin{equation}\label{betatdef}
\frac{\beta(t)}{\alpha\big(\beta(t)\big)^2} \sim d\, \log t\, .
\end{equation}
In other words, $\beta(t)$ is the asymptotic inverse of $t\mapsto t/\alpha(t)^2$ evaluated at
$d \,\log t$, which by \cite[Theorem~1.15.12]{BGT87} exists and is slowly varying.
In order to avoid technical inconveniences, we
assume that the field $\xi$ is bounded from \emph{below}.
See Remark~\ref{essinf} for comments on this issue.

\begin{theorem}[Almost sure asymptotics]\label{asasy}
Suppose Assumptions (H) and (K) hold, and assume that we are in case~(3), i.e.,
$\gamma=1$ and $\kappa^*=0$. Furthermore, suppose that $\beta$ is defined by \eqref{betatdef}
and that $\essinf \xi(0)>-\infty$. Let $\rho >0$ be as in Proposition~\ref{variation}.
Then, almost surely,
\begin{equation}\label{result2}
\frac 1t\log U(t) = \frac{H\big(\beta(t) \alpha(\beta(t))^{-d}\big)}{\beta(t)\,\alpha(\beta(t))^{-d}} -
 \frac{1}{\alpha(\beta(t))^2} \,  \big(\rho (d -\sfrac d2 \log \sfrac{\rho}{\pi } +
\log\sfrac{\rho}{e})+o(1)\big),
\quad \mbox{ as } t\uparrow\infty.
\end{equation}
\end{theorem}

\begin{remark}[The constant]
In Section~\ref{Varprob}, we will see that also the constant
$\rho (d -\sfrac d2 \log \sfrac{\rho}{\pi } +
\log\sfrac{\rho}{e})$ arises as a variational
problem. A remarkable fact is that the first two leading contributions
to $U(t)$ are {\it deterministic}.
\hfill$\Diamond$
\end{remark}

\begin{remark}[Interpretation]
Heuristically, $\alpha(\beta(t))$ is the order of the diameter of the intermittent islands, which
almost surely carry most of the mass of $U(t)$. Note that $\beta(t) = (\log t)^{1+o(1)}$ and
$\alpha(\beta(t)) = (\log t)^{o(1)}$, i.e.,~the size of the intermittent islands increases extremely
slowly. The crucial point in the proof of Theorem~\ref{asasy} is to show the existence of an island
with radius of order $\alpha(\beta(t))$ within the box $[-t,t]^d$ on which the shape of the
vertically shifted and rescaled potential is optimal, i.e., resembles a certain parabola.
To prove this, we use the first moment asymptotics at time $\beta(t)$ locally on that island.
The exponential rate, which is $\beta(t)/\alpha(\beta(t))^2$ has to be balanced against the
\emph{number} of possible islands, which has exponential rate $d \log t$, cf.~\eqref{betatdef}.
\hfill$\Diamond$
\end{remark}

\begin{remark}[Lower tails of the potential]\label{essinf} The assertion of Theorem~\ref{asasy}
remains true {\it mutatis mutandis} if the assumption $\essinf \xi(0)>-\infty$ is replaced, in $d\geq 2$,
by the assumption that $\Prob\{\xi(0)>-\infty\}$ exceeds the critical nearest-neighbour site percolation
threshold. This ensures the existence of an infinite component in the set
${\mathcal C}=\{z\in\Z^d\colon\xi(z)>-\infty\}$, and thus \eqref{result2} holds conditional
on the event that the origin belongs to the infinite cluster in ${\mathcal C}$.
In $d=1$, an infinite cluster exists if and only if $\Prob\{\xi(0)>-\infty\}=1$.
If we assume that $\xi(0)>-\infty$ almost surely and $\langle\log(-\xi(0)\lor 1)\rangle<\infty$,
\eqref{result2} is true verbatim, while otherwise the rate of the almost sure asymptotics
depends on the lower tails of $\xi(0)$; see \cite{BK01a} for details. The effect of
the assumption is to ensure sufficient connectivity in the sense that the mass flow from the origin
to regions where the random potential assumes high values and an approximately optimal shape is not
hampered by deep valleys on the way.

We decided to detail the proof of the almost sure asymptotics under the stronger assertion that
$\essinf \xi(0)>-\infty$. See \cite[Section~5.2]{BK01} for the proof of the analogous assertion in the bounded-potential case under the weaker assumptions. The arguments given there can be extended with
some effort to the situation of the present paper.
\hfill$\Diamond$
\end{remark}

The proof of Theorem~\ref{asasy} is given in Section~\ref{sec-proofquen}.
It essentially follows the strategy of~\cite{BK01}.

\subsubsection{Examples.}\label{sec-ex}
We now explain what kind of upper tail behaviour is covered by
the almost bounded case, arguing separately for the
bounded and unbounded case, denoted by (B) and (U), respectively.
Suppose the distribution of the field $\xi(0)$ satisfies
    \begin{equation}
    \label{example}
    \log \Prob\big\{\xi(0)>r\big\}\sim -e^{f(r)},\qquad
\mbox{as }\begin{cases}r\uparrow\infty&\mbox{in case (U)},\\
    r\uparrow 0={\rm ess sup}\; \xi(0),&\mbox{in case (B)}.
    \end{cases}
    \end{equation}
Here $f$ is a positive, strictly increasing smooth function satisfying
$f'(r)\uparrow \infty$ as $r\uparrow\infty$ in case (U) and $f'(r)r\uparrow \infty$ as
$r\uparrow 0$ in case (B). Note that typical representatives of case~(2)
of the four universality classes are $f(r)\approx cr$  as $r\uparrow\infty$,
violating the condition in case (U); and
typical representatives of case~(4) of the four universality classes
are $f(r)\approx -\sfrac \gamma{1-\gamma} \log|r|$
as $r\uparrow 0$,  violating the condition in case (B). The cumulant generating function behaves like
    \begin{equation}\label{Hasy2}
    H(t)\approx \log\int e^{tr} \exp\bigl\{-e^{f(r)}\bigr\}\, dr
    \approx \sup_r\bigl[tr-e^{f(r)}\bigr]=t r(t)-e^{f(r(t))},
    \end{equation}
where $r(t)$ is asymptotically, as $t\uparrow\infty$, defined via $t=f'(r(t))e^{f(r(t))}$.
Note that $r(t)\uparrow \infty$ in case (U), while $r(t)\uparrow 0$ in case (B),
as $t\uparrow\infty$. Hence, $f'(r(t))\uparrow\infty$ in case (U), while
$f'(r(t))r(t)\uparrow \infty$ in case (B). Rewriting the definition of $r(t)$ as
    $$
    e^{f(r(t))}=\frac{tr(t)}{f'(r(t))r(t)}=o(tr(t)),
    $$
we thus obtain that the first term on the right hand side of \eqref{Hasy2} dominates the second term.
Therefore, we can approximate ${H(t)}/t\approx r(t)$, as $t\uparrow\infty$.
We next assume that $f'(r(\,\cdot\,))$ is slowly varying at
infinity. We then see that, using the fact that $r(t)=f^{-1}\big( \log\frac{t}{f'(r(t))}\big)$
in the last equality,
    $$\begin{aligned}
    H(ty)-yH(t) & \approx ty \Big( f^{-1}\big( \log \sfrac{ty}{f'(r(ty))}\big)-
    f^{-1}\big( \log\sfrac{t}{f'(r(t))}\big) \Big) \\
    &\approx ty \Big( f^{-1}\bigl(\log\sfrac{t}{f'(r(t))}+\log y\bigr)-
    f^{-1}\big(\log\sfrac{t}{f'(r(t))}\big) \Big) \\
    &\approx t\,(y\log y)\, (f^{-1})'\big(\log \sfrac{t}{f'(r(t))}\big)= (y\log y)\, \frac t{f'(r(t))}.
    \end{aligned}$$
Using Proposition~\ref{variation}, this means that the scaling relation in
\eqref{basic} is satisfied with $\kappa(t)=t/{f'(r(t))}$ and $\rho=1$.
As \mbox{$f'(r(t))\uparrow\infty$} is slowly varying,
we see that we are in case~(3) of the four universality classes.

\subsection{Heuristic derivation of Theorem~\ref{momasy}}\label{MA}

\noindent In this section, we give a heuristic explanation  of Theorem~\ref{momasy} in terms of
large deviations for the scaled potential $\xi$. Our proof of Theorem~\ref{momasy}
follows a different strategy.

We use the setup and notation of Section~\ref{sec-ex} and handle the cases (B) respectively (U)
simultaneously. Consequently, the definition \eqref{alphadef} of $\alpha(t)$ reads
    \begin{equation}
    \label{alphadef-examp}
    \alpha(t)^2=\frac {t\alpha(t)^{-d}}{\kappa(t\alpha(t)^{-d})}=f'\bigl(r(t\alpha(t)^{-d})\bigr).
    \end{equation}
We introduce the shifted, scaled potential
    \begin{equation}\label{shirescpot}
\begin{aligned}
    \overline\xi_t(x) & :=\alpha(t)^2\Bigl[\xi\bigl(\lfloor x
    \alpha(t)\rfloor\bigr)-\frac{H(t\alpha(t)^{-d})}{t\alpha(t)^{-d}}\Bigr]\\
    & \approx\alpha(t)^2\Bigl[\xi\bigl(\lfloor x \alpha(t)\rfloor\bigr)-
    r(t\alpha(t)^{-d})+ \sfrac{\alpha(t)^d}{t} e^{f(r(t\alpha(t)^{-d}))} \Bigr],
    \end{aligned}
    \end{equation}
for $x\in Q_{\sR}=[-R,R]^d$.
The process $\overline \xi_t$ satisfies a large deviation principle,
for every $R>0$, on the cube $Q_{\sR}$ with rate
$t\alpha(t)^{-2}$ and rate function $\varphi\mapsto \int_{Q_{\sR}} e^{\varphi(x)-1}\, dx$.
Indeed, with $B_{\sR}=[-R,R]^d\cap \Z^d$,
    $$
    \begin{aligned}
    \Prob\bigl\{\overline\xi_t\approx\varphi\mbox{ on }Q_{\sR}\bigr\}
    &\approx\prod_{z\in B_{R\alpha(t)}}\Prob\Big\{\xi(0)\approx
    \sfrac{\varphi(z\alpha(t)^{-1})}{\alpha(t)^{2}}
    +r(t\alpha(t)^{-d})-\sfrac{\alpha(t)^d}{t} e^{f(r(t\alpha(t)^{-d}))}\Big\}\\
    &\approx \prod_{z\in B_{R\alpha(t)}}\exp\Bigl\{-\exp\Bigl[f\Bigl(r(t\alpha(t)^{-d})
    +\sfrac{\varphi(z\alpha(t)^{-1})}{\alpha(t)^{2}}-\sfrac{\alpha(t)^d}{t}
    e^{f(r(t\alpha(t)^{-d}))}\Bigr)\Bigr]\Bigr\}
    \end{aligned}
    $$
By a Taylor expansion around $r(t\alpha(t)^{-d})$, using that $s=f'(r(s))e^{f(r(s))}$ for
$s=t\alpha(t)^{-d}$ as well as (\ref{alphadef-examp}), we can continue with
 $$\begin{aligned}
   \Prob\bigl\{\overline\xi_t\approx\varphi\mbox{ on }Q_{\sR}\bigr\}
   &\approx \exp\Bigl\{-\alpha(t)^d\int_{Q_{\sR}}\exp\bigl[f(r(t\alpha(t)^{-d}))+
   \frac{\varphi(x)}{\alpha(t)^{2}}\,f'(r(t\alpha(t)^{-d}))-1\bigr]\, dx\Bigr\}\\
    &= \exp\Bigl\{-\frac{t}{f'(r(t\alpha(t)^{-d}))} \, \int_{Q_{\sR}}e^{\varphi(x)-1} \, dx\Bigr\}\\
    &\approx \exp\Bigl\{-\frac{t}{\alpha(t)^2} \, \int_{Q_{\sR}} e^{\varphi(x)-1}\,dx\Bigr\}.
    \end{aligned}
    $$
The asymptotics of $\langle U(t)^p\rangle$ can now be explained as follows. Note that
$U(t)=u(t,0)$, where $u(t,\,\cdot\,)$ is the solution of the parabolic Anderson model
\eqref{PAM} with initial condition $u(0,\,\cdot\,)=1$. We can approximate $u(t,0)$ by
the $w_t(t,0)$ where $(s,z) \mapsto w_t(s,z)$ is the solution to the initial boundary
value problem \eqref{PAM} with zero boundary condition outside the box $B_t$ and
initial condition $w_t(0,\,\cdot\,)=\one_{B_t}$. Let $\lambda^{\rm d}_t(\xi)$ denote
the principal eigenvalue of $\DeltaD+\xi$ in $\ell^2(B_t)$ with zero boundary condition.
Then an eigenfunction expansion shows that
    $$
    U(t)^p= u(t,0)^p \approx w_t(t,0)^p \approx e^{pt \lambda_t^{\rm d}(\xi)}.
    $$
This already explains why the asymptotics of the $p^{\rm th}$ moments
of $U(t)$ are the same as the asymptotics of the moments of $U(pt)$. We proceed by taking $p=1$.
Now the shift invariance and the asymptotic scaling properties of the discrete Laplace operator
yield that
    $$
    \lambda_t^{\rm d}(\xi)=\sfrac{H(t\alpha(t)^{-d})}{t\alpha(t)^{-d}}+
\lambda_t^{\rm d}\bigl(\alpha(t)^{-2}\overline\xi_t(\lfloor\cdot\,\alpha(t)^{-1}\rfloor)\bigr)\approx
    \sfrac{H(t\alpha(t)^{-d})}{t\alpha(t)^{-d}}+\alpha(t)^{-2}\lambda(\overline\xi_t),
    $$
where $\lambda(\psi)$ denotes the principal eigenvalue of $\Delta+\psi$ in
$L^2(Q_{t\alpha(t)^{-d}})$,
with zero boundary condition. Hence,
    \begin{equation}
    \label{asympheur1}
    \langle U(t)\rangle  \approx e^{H(t\alpha(t)^{-d})\alpha(t)^{d}}\Bigl\langle \exp\Bigl\{\frac
    t{\alpha(t)^2}\lambda(\overline\xi_t)\Bigr\}\Bigr\rangle.
    \end{equation}
Using the large deviation principle for $\overline\xi_t$ with $R=t \alpha(t)^{-d}$, and anticipating that
$\psi\mapsto\lambda(\psi)$ has the appropriate continuity and boundedness properties,
we may use Varadhan's lemma to deduce that
$$
\frac 1t\log \langle U(t)\rangle \approx
\sfrac{H(t\alpha(t)^{-d})}{t\alpha(t)^{-d}}-\frac 1{\alpha(t)^2} \chi,
$$
where $\chi$ is given by
    \begin{equation}
    \label{chidefheur}
    \chi=\inf_\psi\Big\{\int_{\R^d} e^{\psi(x)-1}\, dx-\lambda(\psi)\Big\}.
    \end{equation}
We show in Section~\ref{Varprob} that $\chi$ is equal to  $\rho d(1-\frac 12 \log \sfrac \rho \pi)$.
This completes the heuristic derivation of Theorem \ref{momasy}.
The interpretation of the above heuristics is that the moments of the
total mass $U(t)$ are mainly governed by potentials $\xi$ whose shape is approximately given as
$$\xi(\cdot)\approx \sfrac{H(t\alpha(t)^{-d})}{t\alpha(t)^{-d}}
+\alpha(t)^{-2} \psi\bigl(\,\cdot\,\alpha(t)^{-1}\bigr)$$
where $\psi$ is a minimiser of the formula in \eqref{chidefheur}.

\subsection{Variational representations of the constants in Theorem~\ref{momasy} and~\ref{asasy}}
\label{Varprob}

\subsubsection{The constant in Theorem~\ref{momasy}}

Fix $\rho>0$ and define $\chi(\rho)\in\R$ by
    \begin{equation}\label{chidef}
    \chi(\rho)= \inf_{\heap{g\in H^1(\R^d)}{\|g\|_2=1}}\Big\{ \|\nabla g\|_2^2-\Hcal(g^2)\Big\},
    \end{equation}
where $H^1(\R^d)$ is the usual Sobolev space, $\nabla$ the usual (distributional) gradient, and
    \begin{equation}\label{calHdef}
    \skrih(g^2)=\rho\int_{\R^d} g^2(x)\log g^2(x)\, dx.
    \end{equation}
By the logarithmic Sobolev inequality in (\ref{logSobineq}) below,
$\skrih(g^2)\in[-\infty,\infty)$ is well-defined for $g\in H^1(\R^d)$.
Furthermore, we introduce the Legendre transform of $\skrih$ on $L^2(\R^d)$ and the top
of the spectrum of the operator
$\Delta+\psi$ in~$H^1(\R^d)$,
\begin{equation}\label{Legendre}
{\mathcal L}(\psi)=\sup_{g\in L^2(\R^d)}\big\{\langle g^2,\psi \rangle-{\mathcal
H}(g^2)\big\}\qquad\mbox{and}\qquad
\lambda(\psi)=\sup_{\heap{g\in H^1(\R^d)}{\|g\|_2 =1}}\big\{ \langle\psi, g^2\rangle-\|\nabla g\|_2^2\big\}.
\end{equation}
Introduce the functions
    \begin{equation}\label{Gauss}
    g_\rho(x)=\Big(\frac\rho\pi\Big)^{\frac d4}
    e^{-\frac\rho2 |x|^2}\qquad\mbox{and}\qquad\psi_\rho(x)=\rho+
    \rho\frac d2\log\frac\rho\pi-\rho^2|x|^2,\qquad \mbox{for }x\in\R^d.
    \end{equation}
Note that the Gaussian density $g_\rho$ is the unique $L^2$-normalized positive eigenfunction of the
operator $\Delta+\psi_\rho$ in $H^1(\R^d)$ with eigenvalue $\lambda(\psi_\rho)=\rho-\rho d+\rho\frac d2
\log\frac\rho\pi$. It satisfies $\Lcal(\psi_\rho)=\rho$.

\begin{prop}[Solution of the variational formula in \eqref{chidef}]\label{lem-chiident}
For any $\rho\in(0,\infty)$, the infimum in \eqref{chidef} is, up to horizontal shift,
uniquely attained at $g_\rho$.
In particular, $\chi(\rho)=\rho d\big(1-\frac 12 \log \sfrac \rho \pi\big)$
is the constant appearing in Theorem \ref{momasy}. Moreover, $\Lcal$ is identified as
    \begin{equation}\label{Lcalident}
    {\mathcal L}(\psi)=\frac{\rho}e\int_{\R^d} e^{\frac 1\rho \psi(x)}\, dx,
    \end{equation}
and the
\lq dual\rq\ representation is
    \begin{equation}\label{chiident}
    \chi(\rho)=\inf_{\psi\in \skric(\R^d)} \big\{ {\mathcal L}(\psi)-\lambda(\psi)\big\},
    \end{equation}
where $\skric(\R^d)$ is the set of continuous functions $\R^d\to\R$.
Up to horizontal shift, the infimum in \eqref{chiident} is uniquely attained at the
parabola $\psi_\rho$ in \eqref{Gauss}.
\end{prop}

\begin{Proof}{Proof.}
By the {\it logarithmic Sobolev inequality\/} in the form of  \cite[Th.~8.14]{LL01} with
$a=\sqrt{\pi/\rho}$, we have
    \begin{equation}
    \label{logSobineq}
    \|\nabla g\|_2^2 \geq \rho\int_{\R^d} g^2(x)\log g^2(x)\, dx
    + \rho d\big(1-\sfrac 12 \log \sfrac \rho \pi\big),
    \end{equation}
with equality exactly for the Gaussian density $g_\rho$ and its horizontal
shifts. This proves the first statement.
In order to see that \eqref{Lcalident} holds, use Jensen's inequality for any $g\in L^2(\R^d)$ to obtain
\begin{equation}\label{JensenLcal}
\langle g^2,\psi\rangle -\Hcal(g^2)=\rho\|g\|_2^2\int\frac{g^2}{\|g\|_2^2}
\log\frac{e^{\frac 1\rho \psi}}{g^2}\leq \rho\|g\|_2^2\log \frac{\int e^{\frac 1\rho \psi}}{\|g\|_2^2}.
\end{equation}
Equality holds if and only if $g^2=C e^{\frac 1\rho \psi}$ for some $C>0$.
The right side of \eqref{JensenLcal} is maximal precisely for  $\|g\|_2^2=\frac 1e \int
e^{\frac 1\rho \psi}$. Substituting this value, we arrive at \eqref{Lcalident}.

To see the last two statements, we use \eqref{Lcalident} and the formula in \eqref{Legendre} for
$\lambda(\psi)$ to obtain, for any $\psi\in \skric(\R^d)$,
    \begin{equation}\label{chiidentproof}
    {\mathcal L}(\psi)-\lambda(\psi)=
    \inf_{\heap{g\in H^1(\R^d)}{\|g\|_2 =1}}\Big(\|\nabla g\|_2^2-\Hcal(g^2)-\rho\int
    g^2\Big[\frac \psi\rho-\log g^2-e^{\frac1\rho \psi-\log g^2-1}\Big]\Big).
    \end{equation}
The term in square brackets is equal to $\theta -e^{\theta-1}$ for $\theta=\frac \psi\rho-\log g^2$.
Since this is nonpositive and is zero only for $\theta=1$, we have that \lq$\leq$\rq\ holds
in \eqref{chiident}. Furthermore, by restricting the infimum over $g$ to strictly positive
continuous functions and interchanging the order of the infima,
we see that
    $$
    \begin{aligned}
    \inf_{\psi\in \skric(\R^d)} \big\{ {\mathcal L}(\psi)-\lambda(\psi)\big\}
    &\leq \inf_{\heap{g\in H^1(\R^d)}{\|g\|_2 =1, g>0}}
    \inf_{\psi\in \skric(\R^d)}\Big(\|\nabla g\|_2^2-\Hcal(g^2)-\rho\int
    g^2\Big[\frac \psi\rho-\log g^2-e^{\frac1\rho \psi-\log g^2-1}\Big]\Big)\\
    &\leq\inf_{\heap{g\in H^1(\R^d)}{\|g\|_2 =1, g>0}}
    \|\nabla g\|_2^2-\Hcal(g^2)=\chi(\rho),
    \end{aligned}
    $$
by substituting $\psi=\rho+\rho\log g^2$, and we use that the maximizer $g$ of the right hand side
is strictly positive. Therefore, equality holds in \eqref{chiident}. We also know that, by uniqueness of the
solution in \eqref{chidef}, the unique minimizer in \eqref{chiident} is $\psi=\rho+\rho\log g_\rho^2
=\psi_\rho$.
\qed
\end{Proof}

\begin{remark}[Interpretation] Both representations \eqref{chidef} and \eqref{chiident} may be
interpreted in terms of optimal rescaled profiles for the moment asymptotics of the total mass $U(t)$.
While the minimizer $\psi_\rho$ in \eqref{chiident} describes the shape of the potential $\xi$
(see Section~\ref{MA}), the minimizer $g_\rho$ in \eqref{chidef} describes the solution $u(t,\cdot)$,
cf.~Remark~\ref{rem-chiinter}.
\hfill$\Diamond$
\end{remark}

\subsubsection{The constant in Theorem~\ref{asasy}}\label{asvarform}

\noindent We now turn to the variational representation of the constant
appearing in Theorem~\ref{asasy}. We define $\widetilde \chi(\rho)$ by
    \begin{equation}\label{chitildedef}
    \widetilde \chi(\rho)=\inf\{-\lambda(\psi)\colon \psi\in \skric(\R^d),\skril(\psi)\leq 1\}.
    \end{equation}
where we recall that $\Ccal(\R^d)$ is the set of continuous functions $\R^d\to\R$.

\begin{prop}[Solution of the variational formula in \eqref{chitildedef}]
\label{asconstant}
For any $ \rho\in(0,\infty)$, the function $\psi_\rho-\rho \log \frac{\rho}{e}$,
with $\psi_\rho$ as defined in \eqref{Gauss},
is the unique minimizer in \eqref{chitildedef}, and
$\widetilde\chi(\rho )=\chi(\rho ) + \rho \log \sfrac{\rho}{e}$.
\end{prop}

\begin{Proof}{Proof.} Obviously, the condition $\skril(\psi)\leq 1$ in \eqref{chitildedef} may be
replaced by $\skril(\psi)= 1$. In the representation
$$
\widetilde \chi(\rho)=\inf\Big\{\rho\log \Lcal(\psi)-\lambda(\psi)\colon
\psi\in \skric(\R^d),\skril(\psi)= 1\Big\}
$$
we may omit the condition $\skril(\psi)= 1$ completely
since $\rho\log \Lcal(\psi)-\lambda(\psi)$ is invariant under
adding constants to~$\psi$. We use the definition of $\lambda(\psi)$
in \eqref{Legendre} and obtain, after interchanging the infima,
\begin{equation}
\label{tildechirewr}
\widetilde \chi(\rho)=\inf_{\heap{g\in H^1(\R^d)}{\|g\|_2=1}}
\Big\{\|\nabla g\|_2^2-\sup_{\psi\in\Ccal(\R^d)}
\Big(\langle \psi, g^2\rangle-\rho\log \int e^{\frac 1\rho \psi(x)}\, dx\Big)\Big\}
+\rho\log\frac \rho {e}.
\end{equation}
The supremum over $\psi$ is uniquely (up to additive constants) attained at $\psi=\rho\log g^2$ with value $\Hcal(g^2)$, as an application of Jensen's inequality shows:
$$
\begin{aligned}
\rho \log \int e^{\frac 1\rho \psi(x)}\, dx&=\rho \log \int dx\, g^2(x)\,e^{\frac 1\rho \psi(x)-\log g^2(x)}\geq \rho \int dx\, g^2(x)\,\Big(\frac 1\rho \psi(x)-\log g^2(x)\Big)\\
&=\langle \psi, g^2\rangle -\Hcal(g^2).
\end{aligned}
$$
Hence, $\widetilde \chi (\rho)=\chi(\rho)+ \rho \log \frac{\rho}{e}$.
Since $g_\rho$ is, up to horizontal shifts, the unique minimiser in \eqref{chidef},
$\widetilde \psi_\rho=\rho \log g_\rho^2 + C$ is the unique minimizer in (\ref{tildechirewr}).
By the above reasoning, $\widetilde \psi_\rho$ is the unique
minimizer of \eqref{chitildedef}, where $C=-\rho \log \sfrac{\rho}{e}$ is determined
by requiring that ${\mathcal L}(\widetilde \psi_\rho)=1$.
\qed\end{Proof}

\begin{remark}[Interpretation] There is an interpretation of the minimiser of \eqref{chitildedef} in
terms of the optimal rescaled profile of the potential $\xi$ for the almost-sure asymptotics of the
total mass $U(t)$. Indeed, the condition $\Lcal(\psi)\leq 1$ guarantees that, almost surely for all
large $t$,  the profile $\psi$ appears in some \lq microbox\rq\ in the rescaled landscape $\xi$ within
the \lq macrobox\rq\ $B_t=[-t,t]^d\cap\Z^d$, which is one of the intermittent islands. The logarithmic
rate of the total mass, $\frac 1t \log U(t)\approx \lambda_{B_t}(\xi)$, can be bounded from below against
the eigenvalue of $\xi$ in the microbox, which is described by $\lambda(\psi)$. Optimising over all
admissible $\psi$ explains the lower bound in \eqref{result2}. Our proof of the lower bound
in Section~\ref{sec-proofquen} makes this heuristics precise.

The Gaussian density $g_\rho$ in \eqref{Gauss} is the unique positive $L^2$-normalized eigenfunction
of $\Delta+\psi_\rho-\rho \log \sfrac{\rho}e$ corresponding to the eigenvalue $-\widetilde{\chi}(\rho)=
\lambda(\psi_\rho-\rho \log \sfrac{\rho}e)$. It describes the rescaled shape of the solution $u(t,\cdot)$
in the intermittent island. An interesting consequence is that the appropriately rescaled potential and
solution shapes are identical for the moment asymptotics and for the almost sure asymptotics. This
phenomenon also occurs in the cases of the double-exponential distribution and the potentials bounded
from above.
\hfill$\Diamond$
\end{remark}

\subsubsection{Convergence of the variational problem in \eqref{Bikoe}}
We close this section by showing that the variational problem in \eqref{Bikoe}
converges to the variational problem in \eqref{chicase3} as $\gamma\uparrow 1$.
We define
\begin{equation}\label{Bikoe2}
\chi(\rho,\gamma)= \inf_{\heap{g\in H^1(\R^d)}{\|g\|_2=1}} \Big\{ \int_{\R^d}|\nabla g(x)|^2\, dx +
\rho\int_{\R^d}\frac {g^{2\gamma}(x)-g^2(x)}{1-\gamma}\, dx \Big\},
\end{equation}
which is equal to the variational problem in \eqref{Bikoe}.

\begin{prop}[Convergence of the variational problem in \eqref{Bikoe2}]
\label{prop-convgammaVP}
For any $ \rho\in(0,\infty)$,
\begin{equation}
\label{convgammaVP}
\lim_{\gamma\uparrow 1} \chi(\rho,\gamma)=\chi(\rho).
\end{equation}
\end{prop}

\begin{Proof}{Proof.} The upper bound in \eqref{convgammaVP} follows by
substituting the Gaussian density $g=g_\rho$ in \eqref{Gauss} into
the infimum in \eqref{Bikoe2}, and by noting that
$$    \lim_{\gamma\uparrow 1}
    \int_{\R^d}\frac {g^{2\gamma}_\rho(x)-g^2_\rho(x)}{\gamma-1}\, dx
    =\int_{\R^d}g^2_\rho(x)\log{g^2_\rho(x)} \, dx,$$
by an explicit computation of the integrals involved.

For the lower bound in \eqref{convgammaVP}, we bound, for any $\gamma\in [0,1)$
and $g\in H^1(\R^d)$,
$$  \begin{aligned}
    \int_{\R^d}\frac {g^{2\gamma}(x)-g^2(x)}{1-\gamma}\, dx
    &=\int_{\R^d}g^2(x) \frac {e^{(\gamma-1)\log{g^2(x)}}-1}{1-\gamma}\, dx
    \geq -\int_{\R^d}g^2(x)\log{g^2(x)}\, dx,
    \end{aligned}$$
since $e^{\theta}-1\geq \theta$ for every $\theta\in \R$. Therefore,
$\chi(\rho,\gamma)\leq \chi(\rho)$ for every $\gamma\in [0,1)$.
\qed\end{Proof}

The remainder of the paper is as follows. In Section~\ref{sec-cutting} we present an important
auxiliary result on self-intersections of random walks, which will be used in the proof of
Theorem~\ref{momasy} in Section~\ref{sec-proofann}. The proof of Theorem~\ref{asasy} is given
in Section~\ref{sec-proofquen}. Finally, in Section~\ref{sec-erratum} we use the opportunity to
correct an error in the proof of the moment asymptotics in case~(4) from \cite{BK01}.

\section{An auxiliary result on self-intersections of random walks}\label{sec-cutting}

In this section we provide a result on $q$-fold self-intersections of random walks,
for small $q>1$, which is an important tool in the proof of the upper bound in Theorem~\ref{momasy}.

It is also of independent interest. Let $\ell_t(z)=\int_0^t\delta_z(X(s))\, ds$ denote the local time
at $z$ of the simple random walk $(X(s) \colon s\in[0,t])$ on $\Z^d$ with generator
$\DeltaD$, starting at the origin.

\begin{prop}\label{cutting} Fix $q>1$ such that $q(d-2)<d$ and $R>0$.
Let $\alpha(t)\to\infty$ such that $\alpha(t)=\Ocal(t^{2/(2d+2)-\eps})$ for some $\eps>0$. Then
\begin{equation}\label{pnormbound}
\limsup_{\theta\downarrow 0}\limsup_{t\uparrow\infty}\frac{\alpha(t)^2}t\log
\E\Bigl[\exp\Big\{ \theta\alpha(t)^{-\frac 1q [d+(2-d)q]} \Bigl(\sum_{x\in\Z^d}
\ell_t(x)^{q}\Bigr)^{\frac 1{q}}\Big\}
\one\{\supp(\ell_t)\subseteq B_{R\alpha(t)}\}\Bigr]=0.
\end{equation}
\end{prop}

\begin{remark} The result is better understood when rephrasing it in terms of the normalised and
rescaled local times, $L_t(\cdot)=\frac 1t\alpha(t)^d\ell_t(\lfloor\,\cdot\,\alpha(t)\rfloor)$.
Then the exponent may be rewritten as
$$
\alpha(t)^{-\frac 1q [d+(2-d)q]} \Bigl(\sum_{x\in\Z^d}
\ell_t(x)^{q}\Bigr)^{\frac 1{q}}=\frac t{\alpha(t)^2}\|L_t\|_q,
$$
where $\|\cdot\|_q$ is the norm on $L^q(\R^d)$. Hence, \eqref{pnormbound} is a large deviations
result for the $q$-norm of $L_t$ on the scale $t/\alpha(t)^2$. It is known that $(L_t \colon t>0)$
satisfies a large deviation principle on this scale in the weak topology generated by bounded
continuous functions, see for example \cite{GKS04}. However, \eqref{pnormbound} does not follow
from a routine application of Varadhan's lemma, since the $q$-norm is neither bounded nor continuous
in this topology. See \cite{Ch04} for an analogous result for a smoothed version of $L_t$.
\hfill$\Diamond$
\end{remark}

\begin{remark} Our proof yields~\eqref{pnormbound} also without indicator on
$\{\supp(\ell_t)\subseteq B_{R\alpha(t)}\}$ if the sum is restricted to a finite subset of $\Z^d$.
It can easily be extended to a large class of random walks, also in discrete time.
The proof is based on a combinatorial analysis of the
high integer moments of the random variable $\sum_x \ell_t(x)^q$. This method is of crucial
importance in the analysis of intersections and self-intersections of random paths 
\cite{KM02}, and of random walk in random scenery~\cite{GKS04}.
\hfill$\Diamond$
\end{remark}

\begin{Proof}{Proof of Proposition \ref{cutting}.} By $B$ we denote the box
$B=B_{R\alpha(t)}=[-R\alpha(t),R\alpha(t)]^d\cap\Z^d$. In the exponent on the left side
of \eqref{pnormbound}, we restrict the sum to $x\in B$ and forget about the indicator on
$\{\supp(\ell_t)\subseteq B_{R\alpha(t)}\}$. In the following we write $\|\,\cdot\,\|_q$
for the norm in $\ell^q(B)$.

In a first step we reduce the problem to a problem on asymptotics of high integer moments.
Suppose first that there are constants $T,C>0$ such that
    \begin{equation}\label{momentasy}
    \E\bigl[\|\ell_t\|_q^{kq}\bigr]\leq k^{kq} C^{kq}\alpha(t)^{k[d+(2-d)q]},
    \qquad\mbox{for any }t\geq T, k\geq\frac t{\alpha(t)^2}.
    \end{equation}
We now show that this assumption implies \eqref{pnormbound}.
Expanding the exponential series, we rewrite
    $$
    \begin{aligned}
    \E\Bigl[\exp\Big\{\theta & \alpha(t)^{-\frac 1q [d+(2-d)q]} \|\ell_t\|_q\Big\}\Bigr]
    & =\sum_{k=0}^{\infty}\frac{1}{k!} \Big(\theta \alpha(t)^{-\frac 1q[d+(2-d)q]}\Big)^k
    \E\big[\|\ell_t\|_q^{k}\big].
    \end{aligned}
    $$
Abbreviate $k_t=q\lceil t/\alpha(t)^2 \rceil$. Under our assumption,
    \begin{equation}\label{momentasy2}
    \E\bigl[\|\ell_t\|_q^{k}\bigr]
\leq \Big(\frac{k}{q}\Big)^{k} C^{k}\alpha(t)^{\frac{k}{q}[d+(2-d)q]},\qquad \mbox{ for }t\geq T, k\geq k_t,
    \end{equation}
and hence we obtain
    \begin{equation}\label{momentasy2a}
    \begin{aligned}
    \E\Bigl[&\exp\Big\{\theta\alpha(t)^{-\frac 1q [d+(2-d)q]} \|\ell_t\|_q\Big\}\Bigr]
&\leq\sum_{k=0}^{k_t-1} \frac{1}{k!} \theta^k \alpha(t)^{-\frac kq[d+(2-d)q]}
    \E\big[\|\ell_t\|_q^{k}\big] +\sum_{k=k_t}^\infty \frac{1}{k!} \Big(\frac{\theta C k}q\Big)^k.
    \end{aligned}
    \end{equation}
For all sufficiently small $\theta>0$, the second term is estimated as follows:
$$
\sum_{k=k_t}^\infty\frac{1}{k!} \Big(\frac{\theta C k}q\Big)^k\leq
\sum_{k=k_t}^\infty\Big(\frac{\theta C }{eq}\Big)^k
=\frac{\Big(\frac{\theta C}{eq}\Big)^{k_t}}{1-\frac{\theta C }{eq}},
$$
and the exponential rate (in $t\alpha(t)^{-2}$) of the right hand side tends to
$-\infty$ as $\theta\downarrow0$.

For the first term, we bound, using H\"older's inequality and \eqref{momentasy2}, for $k\leq k_t$,
$$    \begin{aligned}
    \E\big[\|\ell_t\|_q^{k}\big]
    &\leq \E\big[\|\ell_t\|_q^{k_t}\big]^{\frac{k}{k_t}}
\leq \Big(\Big(\frac{k_t}{q}\Big)^{k_t} C^{k_t}\alpha(t)^{\frac{k_t}{q}[d+(2-d)q]}\Big)^{\frac{k}{k_t}}
   & =\Big(\frac{k_t}{q}\Big)^{k} C^{k}\alpha(t)^{\frac{k}{q}[d+(2-d)q]}.
    \end{aligned}$$
Therefore, the first term in \eqref{momentasy2a} is bounded by
    $$    \sum_{k=0}^{k_t} \frac{1}{k!} \Big(\frac{\theta C k_t}{q}\Big)^k
    \leq e^{\frac{\theta C}{q} k_t}.$$
This proves that \eqref{momentasy} implies the statement \eqref{pnormbound}.
Therefore, it suffices to prove
\eqref{momentasy} with some constants $C,T>0$. We
use $C$ to denote a generic constant which depends on $R$, $d$ and $q$,
but not on $k$ and $t$, and $C$ may change its value from appearance to
appearance.

To prove \eqref{momentasy}, we write $A_k$ for the set
of maps $\beta\colon B\to\N_0$ satisfying
$\sum_{x\in B}\beta_x=k$. First we write out
    \begin{equation}\label{startcut}
    \begin{aligned}
    \E\bigl[\|\ell_t\|_q^{kq}\bigr]
    &=\sum_{z_1,\dots,z_k\in B}\E\Bigl[\prod_{x\in B}\ell_t(x)^{q\#\{i\colon z_i=x\}}\Bigr]\\
    &=\sum_{\beta\in A_k}\E\Bigl[\prod_{x\in B}\ell_t(x)^{q\beta_x}\Bigr]
    \#\{z\in B^k\colon \beta_x=\#\{z_i=x\}\forall x\}\\
    &=k!\sum_{\beta\in A_k}\E\Bigl[\prod_{x\in B}\ell_t(x)^{q\beta_x}\Bigr]\prod_{x\in B}\frac 1{\beta_x!}.
    \end{aligned}
    \end{equation}
Note that, for $\beta\in A_k$, the numbers $q\beta_x$ are not necessarily integers.
We resolve this problem, in an upper bound, by introducing a further sum over
the set $A_k(\beta)$ of all $\widetilde\beta
\colon B\to\N_0$ satisfying $|\widetilde \beta_x-q\beta_x|<1$ for every $x\in B$.
Then, clearly,
    \begin{equation}\label{startcut1}
    \E\Bigl[\prod_{x\in B}\ell_t(x)^{q\beta_x}\Bigr]\leq \sum_{\widetilde \beta\in A_k(\beta)}
    \E\Bigl[\prod_{x\in B}\ell_t(x)^{\widetilde\beta_x}\Bigr].
    \end{equation}
We fix $\beta\in A_k$ and $\widetilde \beta\in A_k(\beta)$ and denote
$\widetilde k=\sum_{x\in B}\widetilde\beta_x$. Writing out the local times,
we have
    $$
    \E\Bigl[\prod_{x\in B}\ell_t(x)^{\widetilde\beta_x}\Bigr]=\Bigl[\prod_{x\in B}
    \prod_{i=1}^{\widetilde\beta_x}\int_0^td s_i^x\Bigr]\,\P\bigl\{X(s_i^x)=x
    \,\forall x\in B\,\forall i=1,\dots,\widetilde\beta_x\bigr\}.
    $$
The next step is to give new names to the integration variables $s_i^x$ such that
we can order the time variables.  Fix some function
$\varrho\colon \{1,\dots,\widetilde k\}\to B$ such that
$|\varrho^{-1}(\{x\})|=\widetilde\beta_x$ for any $x\in B$.
We continue with, denoting the set of permutations of $1,\dots,\widetilde k$
by $\Sym_{\widetilde k}$,
    \begin{equation}\label{startcut2}
    \begin{aligned}
    \E\Bigl[\prod_{x\in B}\ell_t(x)^{\widetilde\beta_x}\Bigr]
    &=\int_{[0,t]^{\widetilde k}}dt_1\dots dt_{\widetilde k}\,
    \P\{X(t_i)=\varrho(i)\forall i=1,\dots,\widetilde k\}\\
    &=\int_{0<t_1<\dots<t_{\widetilde k}\leq t}dt_1\dots dt_{\widetilde k}
    \sum_{\sigma\in\Sym_{\widetilde k}}\P\bigl\{X(t_{\sigma(i)})=\varrho(i) \;\;\forall i\bigr\}\\
    &= \sum_{\sigma\in\Sym_{\widetilde k}}
    \int_{(0,\infty)^{\widetilde k}}d s_1\dots ds_{\widetilde k}
    \1\Bigl\{\sum_{i=1}^{\widetilde k}s_i\leq t\Bigr\}
    \,\prod_{i=1}^{\widetilde k}p_{s_i}\big(\varrho(\sigma(i-1)),\varrho(\sigma(i))\big),
    \end{aligned}
    \end{equation}
where we switched from $\sigma$ to $\sigma^{-1}$ and substituted $s_i=t_i-t_{i-1}$ (with $t_0=0$),
and we introduced the transition probabilities of continuous time simple random walk,
$p_s(x,y)=\P_x\{X(s)=y\}$. Here we use the convention $\sigma(0)=0$ and $\varrho(0)=0$,
the starting point of the random walk.

We estimate the indicator on the right hand side of \eqref{startcut2} against
$e^{\lambda t}\prod_{i=1}^{\widetilde k}e^{-\lambda s_i}$
for $\lambda=\alpha(t)^{-2}$. Then we integrate out
over all the $s_i$, to obtain
    \begin{equation}\label{prodform}
    \E\Bigl[\prod_{x\in B}\ell_t(x)^{\widetilde\beta_x}\Bigr]\leq e^{\lambda t}
    \sum_{\sigma\in\Sym_{\widetilde k}}\prod_{i=1}^{\widetilde k}
    G_\lambda\bigl(\varrho(\sigma(i-1)),\varrho(\sigma(i))\bigr),
    \end{equation}
where $G_\lambda$ is the Green's function of the walk given by
$G_\lambda(x,y)=\int_0^{\infty} e^{-\lambda s} p_s(x,y)\, ds.$
It will be convenient to use a closed loop of sites, i.e., to
change the convention $\sigma(0)=0$ to the convention
$\sigma(0)=\sigma(\widetilde k)$. This change of conventions
leads to a factor
    \[
    \frac{G_\lambda\bigl(\varrho(0),\varrho(\sigma(1))\bigr)}
    {G_\lambda\bigl(\varrho(\sigma(\widetilde k)),\varrho(\sigma(1))\bigr)},
    \]
which can be bounded by $e^{o(k)}$ since $\sup_{x, y \in B} {G_\lambda(0,y)}/{G_\lambda(x,y)}
\leq e^{o(k)}$, where we recall that $\lambda=\alpha(t)^{-2}$,
$k>t\alpha(t)^{-2}$ and $B=B_{R\alpha(t)}$.

We denote by $P(\widetilde\beta)$ the set of maps $\gamma\colon B\times B\to\N_0$ such that
$\sum_{y\in B}\gamma_{x,y}=\sum_{y\in B}\gamma_{y,x}=\widetilde\beta_x$
for any $x\in B$. Then we can rewrite
    $$\begin{aligned}\E\Bigl[\prod_{x\in B}\ell_t(x)^{\widetilde\beta_x}\Bigr]&\leq e^{\lambda t+o(k)}
    \sum_{\gamma\in P(\widetilde\beta)}\prod_{x,y\in B} G_\lambda(x,y)^{\gamma_{x,y}}
    \sum_{\heap{w\in B^{\widetilde k}}{(w_0=w_{\widetilde k})}}\sum_{\sigma\in\Sym_{\widetilde k}} \\
    &\qquad\qquad\times \1\{\varrho\circ\sigma=w\}\1\big\{\gamma_{x,y}=
    \#\{i\colon w_{i-1}=x,w_i=y\}\;\;\forall x,y\in B\big\}.
    \end{aligned}$$
We can evaluate the sums over $w$ and $\sigma$ using elementary
combinatorics. Indeed, note that
    \begin{equation}
    \label{comb1}
    \#\{\sigma\colon \varrho\circ\sigma=w\}=\prod_{x\in B}\widetilde\beta_x!,
    \end{equation}
since, given $w$ and $\varrho$, the left hand side equals
the number of orders in which one can put $\widetilde k$ objects,
of which $\widetilde\beta_x$ for each $x \in B$ are indistinguishable,
into a row such that the same vector of elements arises. Since one can
only permute within those indices which belong to the same class of
indistinguishable objects, we obtain \eqref{comb1}. This performs the sum
over $\sigma$. To perform the sum over $w$ for fixed $\gamma$,
we use \cite[p.17]{dH00}, to obtain
$$    \#\big\{w\colon \gamma_{x,y}=
    \#\{i\colon w_{i-1}=x,w_i=y\}\;\;\forall x,y\in B\big\}
    \leq \widetilde k \, \prod_{x,y\in B} \frac {\widetilde\beta_x!}{\gamma_{x,y}!}.$$
Therefore, we obtain
    \begin{equation}\label{startcut3}
    \begin{aligned}
   \E\Bigl[\prod_{x\in B}\ell_t(x)^{\widetilde\beta_x}\Bigr]
    & \leq e^{\lambda t+o(k)}\sum_{\gamma\in P(\widetilde\beta)}\prod_{x,y\in B}
    \Bigl[G_\lambda(x,y)^{\gamma_{x,y}}\frac {\widetilde\beta_x!}{\gamma_{x,y}!}\Bigr]
    \prod_{x\in B}\widetilde\beta_x!\\
    &\leq e^{\lambda t+o(k)}\sum_{\gamma\in P(\widetilde\beta)}e^{\sum_{x,y\in B}\gamma_{x,y}}
    \prod_{x,y\in B} \Bigl[\frac{G_\lambda(x,y)^q}{\gamma_{x,y}}\Bigr]^{\frac 1q \gamma_{x,y}}
    \prod_{x,y\in B}\gamma_{x,y}^{\frac {1-q}q \gamma_{x,y}}
    \prod_{x\in B}\widetilde\beta_x^{2\widetilde\beta_x},
    \end{aligned}
    \end{equation}
where we use that $n^n e^{-n}\leq n!\leq n^n$.
We next use that, since $|\widetilde \beta_x-q\beta_x|<1$,
$$    \widetilde{k}=
     \sum_{x,y\in B}\gamma_{x,y}=\sum_{x\in B} \widetilde\beta_x \leq q\sum_{x\in B} \beta_x +
    |B|=qk +|B|.$$
By our assumption on the growth of $\alpha(t)$, we have $|B|\leq C\alpha(t)^d\leq o(k)$ and hence
$e^{\sum_{x,y\in B}\gamma_{x,y}}\leq C^k$. Fix $\gamma\in P(\widetilde\beta)$. We use Jensen's
inequality for the logarithm to obtain
    \begin{equation}\label{Jensen1}
    \begin{aligned}
    \prod_{x,y\in B} \Bigl[\frac{G_\lambda(x,y)^q}{\gamma_{x,y}}\Bigr]^{\frac 1q \gamma_{x,y}}
    &=\exp\Bigl\{\frac 1q \sum_{x\in B}\widetilde\beta_x\sum_{y\in B}\frac{\gamma_{x,y}}{\widetilde\beta_x}
    \log \frac{G_\lambda(x,y)^q}{\gamma_{x,y}}\Bigr\}\\
    &\leq \exp\Bigl\{\frac 1q \sum_{x\in B}\widetilde\beta_x
    \log\Bigl(\sum_{y\in B}\frac{G_\lambda(x,y)^q}{\widetilde\beta_x}\Bigr)\Bigr\}.
    \end{aligned}
    \end{equation}
Recall that $\lambda=\alpha(t)^{-2}$. Since $(d-2)q<d$, there is a constant $C$
(only depending on $R$, $d$ and $q$) such that, for any $x\in B$,
    \begin{equation}\label{Ghochp}
    \sum_{y\in B}G_{\alpha(t)^{-2}}(x,y)^q\leq C \alpha(t)^{d+(2-d)q}.
    \end{equation}
This gives that
    \begin{equation}
    \label{Jensenbd}
    \prod_{x,y\in B} \Bigl[\frac{G_\lambda(x,y)^q}{\gamma_{x,y}}\Bigr]^{\frac 1q \gamma_{x,y}}
    \leq C^k \alpha(t)^{[d+(2-d)q]\widetilde k/q}\prod_{x\in B}\beta_x^{-\frac 1q\beta_x}.
    \end{equation}
We substitute \eqref{Jensenbd} into \eqref{startcut3}
and summarise \eqref{startcut}, \eqref{startcut1} and \eqref{startcut3}.
Using that $|\widetilde k -qk|\leq |B|$, we obtain
    \begin{equation}\label{summarize}
    \E\Bigl[\|\ell_t\|_q^{qk}\Bigr]
    \leq \widetilde k^{\widetilde k}C^k\alpha(t)^{[d+(2-d)q](k+\frac 1q|B|)}\sum_{\beta\in A_k}
    \sum_{\widetilde\beta\in A_k(\beta)}\sum_{\gamma\in P(\widetilde\beta)}
    \prod_{x,y\in B}\Bigl[\frac{(\frac 1{\widetilde k}
    \gamma_{x,y})^{\frac {1-q}q}}{(\frac 1{\widetilde k}\widetilde\beta_x)^{\frac 1q -2}
    (\frac 1k \beta_x)^{\beta_x/\widetilde\beta_x}}\Bigr]^{\gamma_{x,y}}.
    \end{equation}
Note that, by our growth assumption on $\alpha(t)$ and since $k\geq t/\alpha(t)^2$,
$$    \widetilde k^{\widetilde k}\leq (qk)^{qk} C^{|B|}k^{|B|}\leq (qk)^{qk}e^{o(k)}.$$
The product is estimated with the help of Jensen's inequality for the
logarithm, together with the fact that $y\mapsto \frac {\gamma_{x,y}}
{\widetilde \beta_x}$ is a probability measure, as follows:
    $$
    \begin{aligned}
    \prod_{x,y\in B}\Bigl[\frac{(\frac 1{\widetilde k}\gamma_{x,y})^{\frac {1-q}q}}
    {(\frac 1{\widetilde k}\widetilde\beta_x)^{\frac 1q -2}
    (\frac 1k \beta_x)^{\beta_x/\widetilde\beta_x}}\Bigr]^{\gamma_{x,y}}
    &=\exp\Bigl\{\frac{q-1}q\sum_{x\in B}\widetilde \beta_x\sum_{y\in B}
    \frac {\gamma_{x,y}}{\widetilde \beta_x}\log\frac{\widetilde\beta_x
    (\frac 1{\widetilde k}\widetilde\beta_y)}{\gamma_{x,y}}\Bigr\}\\
    &\qquad\times\exp\Bigl\{-\sum_{x\in B}\beta_x\log\frac{\beta_x} k
    +\frac 1q\sum_{x\in B}\widetilde \beta_x\log\frac{\widetilde\beta_x}{\widetilde k}\Bigr\}\\
    &\leq \exp\Bigl\{-\sum_{x\in B}\beta_x\log\frac{\beta_x} k
    +\frac 1q\sum_{x\in B}\widetilde \beta_x\log\frac{\widetilde\beta_x}{\widetilde k}\Bigr\}.
    \end{aligned}
    $$
Now recall that $q\beta_x -1\leq \widetilde \beta_x\leq q\beta_x+1\leq 2q\beta_x$
for $\beta_x>0$ to bound
    $$
    \sum_{x\in B}
    \widetilde \beta_x\log\frac{\widetilde\beta_x}{\widetilde k}
    \leq \sum_{x\in B}
    (q\beta_x-1)\log\frac{2q\beta_x}{\widetilde k}
    =q\sum_{x\in B}
    \beta_x\log\frac{\beta_x}{k} +\sum_{x\in B}
    \log\frac{\widetilde k}{2q\beta_x}+qk\log\frac{k}{\widetilde k},
    $$
so that we arrive at
    $$
    \begin{aligned}
    -\sum_{x\in B}\beta_x\log\frac{\beta_x} k+\frac 1q\sum_{x\in B}
    \widetilde \beta_x\log\frac{\widetilde\beta_x}{\widetilde k}&\leq
    k\log\frac{k}{\widetilde k}+\frac 1q \sum_{x\in B}\log\frac{\widetilde k}{2q\beta_x}\\
    &\leq k\log\frac{qk}{\widetilde k} +Ck+C|B|\log k \leq Ck,
    \end{aligned}
    $$
since $qk/\widetilde k$ converges to one and since $|B|\log k\leq C\alpha(t)^d\log k\leq o(k)$.
Hence, we have estimated the product on the right hand side of \eqref{summarize} against $C^k$
uniformly in $\beta\in A_k$, $\widetilde \beta\in A_k(\beta)$ and $\gamma\in P(\widetilde \beta)$.
Our growth condition on $\alpha(t)$ implies that each of the sums can be estimated against
$e^{o(k)}$. Indeed,
    $$
    |P(\widetilde\beta)|\leq k^{|B|^2}\leq e^{C\alpha(t)^{2d}\log k}\leq e^{o(k)},
    $$
for any $\widetilde\beta\in A_k(\beta)$ and for any $\beta\in A_k$. Furthermore,
    $|A_k(\beta)|\leq 2^{|B|}\leq e^{o(k)}$
for any $\beta\in A_k$, and finally $|A_k|\leq k^{|B|}\leq e^{o(k)}$. Therefore, we obtain
$$    \E\bigl[\|\ell_t\|_q^{kq}\bigr]
    \leq C^k k^{qk}\alpha(t)^{[d+(2-d)q]k}\alpha(t)^{C|B|}\leq k^{kq}
    C^{kq}\alpha(t)^{k[d+(2-d)q]},$$
where we again used our growth condition on $\alpha(t)$. This completes the proof.
\qed\end{Proof}

\section{The moment asymptotics: Proof of Theorem~\ref{momasy}}
\label{sec-proofann}

\noindent Our analysis is based on the link between the random-walk and
random-field descriptions provided by the Feynman-Kac formula. Let
$(X(s)\colon s\in[0,\infty))$ be the continuous-time simple random
walk on $\Z^d$ with generator $\DeltaD$. By $\P_z$ and
$\E_z$ we denote the probability measure, respectively, the expectation
with respect to the walk starting at $X(0)=z\in\Z^d$.

For any potential $V\colon \Z^d\to[-\infty,\infty)$, we denote by $u^V$
the unique nonnegative, continuous solution of the initial-value~problem
\begin{equation}\label{AndersonV}
\begin{array}{rcll}
\displaystyle
\partial_t u(t,z)\!\!\! &=&\!\!\!\DeltaD u(t,z)+V(z)
u(t,z),\qquad
&\mbox{ for } (t,z)\in(0,\infty)\times \Z^d,\\
u(0,z)\!\!\!&=&\!\!\!1,&\mbox{ for } z\in\Z^d.
\end{array}
\end{equation}
The Feynman-Kac formula allows us to express $u^V$ as
\begin{equation}\label{FK}
u^V(t,z)=\E_z\Big[\exp\Big\{\int_0^t
V\bigl(X(s)\bigr)\,ds\Big\}\Big],\qquad \mbox{ for }z\in\Z^d,\, t>0.
\end{equation}
Hence, $u^\xi$ is the solution of the parabolic Anderson problem in \eqref{PAM} with initial condition $u(0,z)=1$ for all $z\in\R^d$, and the main object of our study is  $U(t)=u^\xi(t,0)$.
\noindent Introduce the vertically shifted random potential
\begin{equation}\label{vtdef}
\xi_t(z)=\xi(z)- H\Big(\frac{t}{\alpha(t)^d}\Big) \frac{\alpha(t)^d}{t}.
\end{equation}
Note that $t$ is a parameter here, and $\xi_t$  should not be seen as a time-dependent random potential.
Fix $p\in(0,\infty)$. Then Theorem~\ref{momasy} is equivalent to the statement
\begin{equation}\label{momasyaim}
\lim_{t\uparrow\infty}\frac{\alpha(pt)^2}{pt}\log \big\langle u^{\xi_{pt}}(t,0)^p\big\rangle=-\chi(\rho),
\end{equation}
where $\chi(\rho)$ is defined in \eqref{chidef}. We approximate $u^{\xi_{pt}}$ by finite-space versions.
Let $R>0$ and let $B_{\sR}=[-R,R]^d\cap \Z^d$ be the centred box in $\Z^d$ with radius $R$.
Introduce $u^V_{\sR}\colon [0,\infty)\times \Z^d\to[0,\infty)$ by
\begin{equation}
\label{altformula}
u^V_{\sR}(t,z)=\E_z\Bigl[\exp\Big\{\int_0^t
V\bigl(X(s)\bigr)\,ds\Big\}\one\bigl\{\supp(\ell_t)\subseteq
B_{\sR}\bigr\}\Bigr],
\end{equation}
where $\ell_t(z)=\int_0^t\delta_z(X(s))\, ds$ are the local times of the random walk.
Note that $u_r^V\leq u_{\sR}^V\leq u^V$ for $0<r<R<\infty$. In the finite space setting we can
work easily with eigenfunction expansions: We look at the function
    \begin{equation}\label{FKp}
    p_R^V(t,y,z)=\E_y \Bigl[e^{\langle V,\ell_t\rangle}\one\bigl\{\supp (\ell_t)\subseteq
    B_R\bigr\}\one\bigl\{X(t)=z\bigr\}\Bigr]\qquad \mbox{ for }y,z\in\Z^d,
    \end{equation}
and the eigenvalues, $\lambda_1>\lambda_2\geq \lambda_3\geq \dots\geq\lambda_n$, of the operator
$\DeltaD+V$ in $\ell^2(B_R)$ with zero boundary condition, where we abbreviate $n=|B_R|$. We may pick an
orthonormal basis of corresponding eigenfunctions ${\rm e}_k$.
By convention, ${\rm e}_k$ vanishes outside $B_R$. Note that 
$\sum_{z\in B_R}p_R^V(t,y,z)=u_R^V(t,y)$. Furthermore, we have
the eigenfunction expansion
    \begin{equation}\label{Fourierp}
    p_R^V(t,y,z)=\sum_{k} e^{t\lambda_k} {\rm e}_k(y){\rm e}_k(z).
    \end{equation}
In particular,
    \begin{equation}\label{Fourierpu}
    u_R^V(t,z)=\sum_{k} e^{t\lambda_k} \langle{\rm e}_k,\one\rangle{\rm e}_k(z).
    \end{equation}
The following proposition carries out the necessary
large deviations arguments for the case $p=1$, and is the key
result for the proof of \eqref{momasyaim}.

\begin{prop}\label{largedeviat1}\ \\ \vspace{-0.8cm}
\begin{itemize}
\item[(i)] Let $R>0$. Then $\displaystyle\limsup_{t\uparrow\infty}\frac{\alpha(t)^2}t\log\left\langle
 u^{\xi_t}_{R\alpha(t)}(t,0) \right\rangle \leq -\chi(\rho).$
\item[(ii)]$\displaystyle\liminf_{R\uparrow\infty}\liminf_{t\uparrow\infty}
\frac{\alpha(t)^2}t\log\left\langle u^{\xi_t}_{R\alpha(t)}(t,0) \right\rangle \geq -\chi(\rho).$
\end{itemize}
\end{prop}

The proofs of Proposition~\ref{largedeviat1}(i) and (ii)
are deferred to Sections~\ref{sec-pfprops1} and~\ref{sec-pfprops2}, respectively.

\subsection{Proof of {\bf{(\ref{momasyaim})}} subject to Proposition~\ref{largedeviat1}}\label{sec-FKform}

\begin{Proof}{Proof of the lower bound in {\bf{(\ref{momasyaim})}}.}
All we have to do is to show that, as $t\uparrow\infty$,
    \begin{equation}
    \label{start}
    \bigl\langle u^{\xi_{pt}}(t,0)^p\bigr\rangle\geq e^{o(t\alpha(pt)^{-2})}
    \Big\langle u_{R\alpha(pt)}^{\xi_{pt}}(pt,0)\Big\rangle.
    \end{equation}
To prove this, we repeat the proof of \cite[Lemmas~4.1 and 4.3]{BK01} for the reader's convenience. We
abbreviate $r=R\alpha(pt)$, $V=\xi_{pt}$, $u=u^{V}$, $u_r=u^{V}_r$ and $p_r=p_r^V$. Note that
$|B_r|=e^{o(t\alpha(pt)^{-2})}$.

Now we prove \eqref{start}. First we assume that $p\in(0,1)$. Use the shift invariance of the
distribution of the field $V$ and the inequality $\sum_i x_i^p\geq (\sum_i x_i)^p$ for
nonnegative $x_i$ to estimate
    \begin{equation}\label{startproof1}
    \bigl\langle u(t,0)^p\bigr\rangle=\Big\langle \frac 1{|B_r|}\sum_{z\in B_r}u(t,z)^p\Big\rangle
    \geq e^{o(t\alpha(pt)^{-2})}\Big\langle\sum_{z\in B_r}u(t,z)^p\Big\rangle
    \geq e^{o(t\alpha(pt)^{-2})}
    \Big\langle\Big(\sum_{z\in B_r}u(t,z)\Big)^p\Big\rangle.
    \end{equation}
By $\|\cdot\|$ we denote the norm on $\ell^2(B_r)$. According to Parseval's identity,
the numbers $\langle {\rm e}_k,\one\rangle^2/\|\one\|^2$ sum up to one. Using the Fourier expansion
in \eqref{Fourierpu} and Jensen's inequality, we obtain
    \begin{equation}\label{startproof2}
    \begin{aligned}
    \Big\langle\Big(\sum_{z\in B_r}u(t,z)\Big)^p\Big\rangle&
    \geq\|\one\|^{2p}\Big\langle\Big(\sum_k e^{t\lambda_k}
    \frac{\langle {\rm e}_k,\one\rangle^2}{\|\one\|^2}\Big)^p\Big\rangle
    \geq e^{o(t\alpha(pt)^{-2})}\Big\langle\sum_k e^{pt\lambda_k}
    \frac{\langle {\rm e}_k,\one\rangle^2}{\|\one\|^2}\Big\rangle\\
    &\geq e^{o(t\alpha(pt)^{-2})}\Big\langle \sum_{z\in B_r}u_r(pt,z)\Big\rangle
    \geq e^{o(t\alpha(pt)^{-2})}\big\langle u_{r}(pt,0)\big\rangle.
    \end{aligned}
    \end{equation}
Substituting \eqref{startproof2} in
\eqref{startproof1} completes the proof of \eqref{start} in the case $p\in(0,1)$.

Now we turn to the case $p\in[1,\infty)$. We use the first equation
in \eqref{startproof1}, Jensen's inequality, the eigenfunction expansion in
\eqref{Fourierpu} and the inequality $(\sum_ix_i)^p\geq \sum_i x_i^p$
for nonnegative $x_i$ to obtain
\begin{equation}\label{startproof3}
\begin{aligned}
\bigl\langle u(t,0)^p\bigr\rangle&\geq\Big\langle \Big(\frac 1{|B_r|}
\sum_{z\in B_r}u(t,z)\Big)^p\Big\rangle
\geq e^{o(t\alpha(pt)^{-2})}\Big\langle\Big(\sum_k e^{t\lambda_k}
\langle {\rm e}_k,\one\rangle^2\Big)^p\Big\rangle\\
&\geq e^{o(t\alpha(pt)^{-2})}\Big\langle\sum_k e^{pt\lambda_k}\langle {\rm e}_k,\one\rangle^{2p}\Big\rangle.
\end{aligned}
\end{equation}
Now we use Jensen's inequality as follows
\begin{equation}\label{startproof4}
\begin{aligned}
\Big\langle\sum_k e^{pt\lambda_k}\langle {\rm e}_k,\one\rangle^{2p}\Big\rangle&=\frac{\big\langle
\sum_k e^{pt\lambda_k}\langle {\rm e}_k,\one\rangle^{2p}\big\rangle}
{\big\langle\sum_k e^{pt\lambda_k}\big\rangle}\,
\Big\langle\sum_k e^{pt\lambda_k}\Big\rangle
\geq \left(\frac{\big\langle\sum_k e^{pt\lambda_k}\langle {\rm e}_k,\one\rangle^{2}\big\rangle}
{\big\langle\sum_k e^{pt\lambda_k}\big\rangle}\right)^p\,\Big\langle\sum_k e^{pt\lambda_k}\Big\rangle\\
&=\Big\langle \sum_{z\in B_r}u_r(pt,z)\Big\rangle\,\left(\frac{\big\langle\sum_k e^{pt\lambda_k}
\big\rangle}{\big\langle\sum_k e^{pt\lambda_k}\langle {\rm e}_k,\one\rangle^{2}\big\rangle}\right)^{1-p}
\geq \big\langle u_{r}(pt,0)\big\rangle.
\end{aligned}
\end{equation}
In the last step, we have used the eigenfunction expansions in
\eqref{Fourierp} and \eqref{Fourierpu} to see that
the ratio is not bigger than one. Combining \eqref{startproof3}
and \eqref{startproof4} completes the proof of \eqref{start}
also in the case $p\in[1,\infty)$.
\qed\end{Proof}

\begin{Proof}{Proof of the upper bound in {\bf{(\ref{momasyaim})}}.}
A main ingredient in our proof is the following preparatory lemma,
which provides, for any potential~$V$,  an estimate of~$u^{V}(t,0)$ in terms of the maximal principal
eigenvalue of $\Delta^{\rm d}+V$ in small subboxes
(`microboxes') of a `macrobox'. For $z\in\Z^d$ and $R>0$, we
denote by $\lambdaD_{z;R}(V)$ the principal eigenvalue of the
operator $\DeltaD+V$ with Dirichlet boundary conditions in
the {\it shifted\/} box $z+B_{R}$.

\begin{lemma}
\label{compact}
Let $r\colon (0,\infty) \to (0,\infty)$ such that $r(t)/t \uparrow \infty$.
For $R,t>0$ let $B_{\sR}(t)=B_{r(t)+2\lfloor R\rfloor}$. Then there
is a constant $C>0$ such that, for any sufficiently large $R, t$ and
any  potential $V\colon\Z^d\to [-\infty,\infty)$,
    \begin{equation}\label{upperbound}
    u^V(t,0)\leq \E\bigl[e^{2\int_0^t V(X_s)\, ds}\bigr]^{1/2}e^{-r(t)}+e^{Ct/R^2}\bigl(3r(t)\bigr)^d
    \exp\left\{t\max_{z\in B_{\sR}(t)}\lambdaD_{z;2R}(V)\right\}.
    \end{equation}
\end{lemma}

\begin{Proof}{Proof.} This is a modification of the proof of \cite[Proposition~4.4]{BK01},
which refers to nonpositive potentials $V$ only. The proof of  \cite[Proposition~4.4]{BK01}
consists of \cite[Lemma~4.5]{BK01} and \cite[Lemma~4.6]{BK01}. The latter states that
\begin{eqnarray}
u^V_{r(t)}(t,0)&\leq &e^{Ct/R^2}\bigl(3r(t)\bigr)^d
    \exp\left\{t\max_{z\in B_{\sR}(t)}\lambdaD_{z;2R}(V)\right\}.
    \label{BK4.6}
\end{eqnarray}
A careful inspection of the proof shows that no use is made
of nonpositivity of $V$ and hence \eqref{BK4.6} applies in the present setting.

In order to estimate $u^V(t,0)-u^V_{r(t)}(t,0)$, we introduce the exit time
$\tau_{\sR}=\inf\{t>0\colon X(t)\notin B_{\sR}\}$ from the box $B_{\sR}$ and use the Cauchy-Schwarz
inequality to obtain
    $$\begin{aligned}
    u^V(t,0)-u^V_{r(t)}(t,0)&=\E\Big[\exp\Big\{\int_0^t V(X(s))\, ds\Big\}\1\{\tau_{r(t)}\leq t\}\Big]\\
    &\leq \E\bigl[e^{2\int_0^t V(X_s)\, ds}\bigr]^{1/2}\P\{\tau_{r(t)}\leq t\}^{1/2}.
    \end{aligned}$$
According to \cite[Lemma~2.5(a)]{GM98}, for any $r>0$,
    $$\P\{\tau_{r}\leq t\}\leq 2^{d+1} \exp\Big\{-r \Big(\log\frac r{dt}-1\Big)\Big\}.$$
Hence, we may estimate $\P\{\tau_{r(t)}\leq t\}^{1/2}\leq e^{-r(t)}$, for sufficiently large $t$,
completing the proof.
\qed
\end{Proof}

We now complete the proof of the upper bound in \eqref{momasyaim}, subject to Proposition~\ref{largedeviat1}.
Let $p\in (0,\infty)$ and fix $R>0$.
First, notice that the second term in \eqref{upperbound}
can be estimated in terms of a sum,
    \begin{equation}
    \label{maxesti}
    \exp\left\{pt\max_{z\in B_{\sR}(t)} \lambdaD_{z;2R}(V)\right\}\leq
    \sum_{z\in B_{\sR}(t)}e^{pt\lambdaD_{z;2R}(V)}.
    \end{equation}
Thus, applying \eqref{upperbound} to $u^{\xi_{pt}}(t,0)$ with $R$ replaced by $R\alpha(pt)$,
raising both sides to the $p$-th power, and using \eqref{maxesti}, we get
$$    \begin{aligned}
    u^{\xi_{pt}}(t,0)^p&\leq 2^p\Bigl(\E\big[e^{2\int_0^t \xi_{pt}(X(s))\, ds}\big]^{p/2} e^{-pr(t)}\\
    &\qquad+  e^{Cpt/(R^2\alpha(pt)^2)}
    \bigl(3r(t)\bigr)^{pd}\!\!\!\sum_{z\in
    B_{R\alpha(pt)}(t)}\!\!\!e^{pt\lambdaD_{z;2R\alpha(pt)}(\xi_{pt})}\Bigr).
    \end{aligned}$$
Next we take the expectation with respect to $\xi$ and note that, by
the shift-invariance of $\xi$, the distribution of $\lambdaD_{z;
2R\alpha(pt)}(\xi)$ does not depend on $z\in\Z^d$. This gives
    \begin{equation}\label{expectation}
    \begin{aligned}
    \langle u^{\xi_{pt}}(t,0)^p\rangle &\leq 2^p\Bigl(\Big\langle
    \E\big[e^{2\int_0^t \xi_{tp}(X(s))\, ds}\big]^{p/2}\Big\rangle e^{-pr(t)}\\
    &\qquad +  e^{Cpt/(R^2\alpha(pt)^2)}
    \bigl(3r(t)\bigr)^{pd+d} \,\big\langle
    e^{pt\lambdaD_{0;2R\alpha(pt)}(\xi_{pt})}\big\rangle\Bigr).
    \end{aligned}
    \end{equation}
In order to show that the first term on the right is negligible, estimate, in the case $p\geq 2$, with the help of Jensen's inequality and Fubini's theorem,
    $$
    \begin{aligned}
    \Big\langle
    \E\big[e^{2\int_0^t \xi_{tp}(X(s))\, ds}\big]^{p/2}\Big\rangle&\leq
    \E\big[\big\langle e^{\frac 1t\int_0^t pt\xi(X(s))\, ds}\big\rangle\big]
    \exp\Big\{-\frac{H(pt\alpha(pt)^{-d})}{\alpha(pt)^{-d}}\Big\}\\
    &\leq \E\Big[\Big\langle\frac 1t\int_0^t e^{ pt\xi(X(s))}\, ds\Big\rangle\Big]
    \exp\Big\{-\frac{H(pt\alpha(pt)^{-d})}{\alpha(pt)^{-d}}\Big\}\\
    &=e^{H(pt)}
    \exp\Big\{-\frac{H(pt\alpha(pt)^{-d})}{\alpha(pt)^{-d}}\Big\}.
    \end{aligned}
    $$
In the case $p<2$, a similar calculation shows that
    $$
    \Big\langle
    \E\big[e^{2\int_0^t \xi_{tp}(X(s))\, ds}\big]^{p/2}\Big\rangle\leq e^{\frac p2 H(2t)}
    \exp\Big\{-\frac{H(pt\alpha(pt)^{-d})}{\alpha(pt)^{-d}}\Big\}.
    $$
Hence, for the choice $r(t)=t^2$, the first term on the right hand side of \eqref{expectation} satisfies
    \begin{equation}\label{firsttermneg}
    \limsup_{t\uparrow \infty}\frac{\alpha(pt)^2}{pt}\log\Big(\Big\langle
    \E\big[e^{2\int_0^t \xi_{tp}(X(s))\, ds}\big]^{p/2}\Big\rangle e^{-pr(t)}\Big)=-\infty,
    \end{equation}
where we use that $H(t)/t$ and $\alpha(t)$ are slowly varying.

In \eqref{expectation}, take the logarithm,
multiply by $\alpha(pt)^2/(pt)$ and let $t\uparrow\infty$. Then we
have that
    \begin{equation}
    \label{momldp}
    \limsup_{t\uparrow\infty}\frac {\alpha(pt)^2}{pt}
    \log\bigl\langle u^{\xi_{pt}}(t,0)^p\bigr\rangle\leq
    \frac C{R^2}+\limsup_{t\uparrow\infty}\frac {\alpha(pt)^2}{pt}\log\bigl\langle
    \exp\{pt\lambdaD_{0;2R\alpha(pt)}(\xi_{pt})\}\bigr\rangle,
    \end{equation}
where we also used that $r(t)^{pd+d}=e^{o(t\alpha(pt)^{-2})}$ as
$t\uparrow\infty$.
Now we estimate the right hand side of \eqref{momldp}. We denote by
$\lambda^{{\rm d},k}_{0;R\alpha(pt)}(\xi_{pt})$ the $k^{\rm th}$
eigenvalue of $\DeltaD + \xi_{pt}$ in the box $B_{R\alpha(pt)}$
with zero boundary condition. Using an eigenfunction expansion
as in \eqref{Fourierp},
we get
 \begin{equation}\label{fromBK}
 \begin{aligned}
 \Big\langle
    \exp\bigl\{pt \lambdaD_{0;R\alpha(pt)}(\xi_{pt})\bigr\}\Big\rangle&\leq \Big\langle\sum_{k}
\exp\bigl\{pt
\lambda^{{\rm d},k}_{0;R\alpha(pt)}(\xi_{pt})\bigr\}\Big\rangle = \sum_{x\in B_{R\alpha(pt)}}
\big\langle p_{R\alpha(pt)}^{\xi_{pt}}(pt,x,x)\big\rangle  \\
&\leq \sum_{x\in B_{R\alpha(pt)}} \big\langle u_{R\alpha(pt)}^{\xi_{pt}}(pt,x)\big\rangle\\
& \leq \sum_{x\in B_{R\alpha(pt)}} \Big\langle\E_x\Big[e^{\int_0^{pt}\xi_{pt}(X(s))\, ds}
\1\{\supp(\ell_{pt})\subseteq x+B_{2R\alpha(pt)}\}\Big]\Big\rangle\\
    &\leq |B_{R\alpha(pt)}|\big\langle u_{2R\alpha(pt)}^{\xi_{pt}}(pt,0)\big\rangle,
    \end{aligned}
    \end{equation}
where we also used the shift-invariance. Recall that $|B_{R\alpha(pt)}|\leq
e^{o(t\alpha(pt)^{-2})}$. We finally use Proposition~\ref{largedeviat1}(i)
for $pt$ instead of $t$ to complete the proof of the upper bound in
\eqref{momasyaim}.
\qed\end{Proof}

\subsection{Proof of Proposition~\ref{largedeviat1}(i)}
\label{sec-pfprops1}

\noindent Recall the local times of the walk, $\ell_t(z)=\int_0^t \one\{X(s)=z\}\,ds$.
Note that $\int_0^t
V(X(s))\,ds=\langle V,\ell_t\rangle$,
where $\langle\,\cdot\,,\,\cdot\,\rangle$ stands for the
inner product on $\ell^2(\Z^d)$.
From \eqref{altformula} with $V=\xi_t$, we have
    \begin{equation}\label{utident1}
    \big\langle u^{\xi_{t}}_{R\alpha(t)}(t, 0) \big\rangle
    =\big\langle\E_0\big[ e^{\langle\xi_t,\ell_t\rangle}\one\{\supp(\ell_t)\subseteq B_{R\alpha(t)}\}\big]\big\rangle.
    \end{equation}
Recall from \eqref{Hdef} that $\langle e^{l\xi(x)}\rangle = e^{H(l)}$ for any
$l\in\R$ and $x\in \Z^d$. We carry out the expectation with respect to the
potential, and obtain, using Fubini's theorem and the independence of the
potential variables,
    \begin{equation}\label{utident2}
    \begin{aligned}
    \big\langle  u^{\xi_{t}}_{R\alpha(t)}(t, 0) \big\rangle
    &= e^{-\alpha(t)^d H(t/\alpha(t)^d)} \me_0
    \Big[ \Big\langle e^{\sum_{x\in R\alpha(t)}\ell_t(x)\xi(x)}\Big\rangle
    \one\{\supp(\ell_t)\subseteq B_{R\alpha(t)}\} \Big] \\
    & = \me_0\Big[ \exp\Big\{ \sum_{x\in B_{R\alpha(t)}}\Big[
    H(\ell_t(x))-\ell_t(x) \sfrac{\alpha(t)^d}{t} H(t/\alpha(t)^d)\Big]\Big\}
    \one\{\supp(\ell_t)\subseteq B_{R\alpha(t)}\} \Big],
    \end{aligned}
    \end{equation}
where we also use that $\sum_{x\in \Z^d}\ell_t(x)=t$. We now split the
sum in the exponent into a part where we have some control over the size
of the local times, and a part with very large local times. Introducing
    \begin{eqnarray}
    \Hcal_{\sM}^{\ssup{t}}(\ell_t)&=&\frac {\alpha(t)^2}t\sum_{x\in B_{R\alpha(t)}}\Big[H(\ell_t(x))-\ell_t(x) \sfrac{\alpha(t)^d}{t} H(t/\alpha(t)^d)\Big]\one\{ \ell_t(x) \leq  \sfrac{Mt}{\alpha(t)^d} \},\\
    \Rcal_{\sM}^{\ssup {t}}( \ell_t)&=&\sum_{x\in B_{R\alpha(t)}}\Big[H(\ell_t(x))-\ell_t(x) \sfrac{\alpha(t)^d}{t} H(t/\alpha(t)^d)\Big]\one\{ \ell_t(x)>  \sfrac{Mt}{\alpha(t)^d} \},\label{remainder}
    \end{eqnarray}
 we have
$$    \big\langle  u^{\xi_{t}}_{R\alpha(t)}(t, 0) \big\rangle
    = \me_0\Big[ \exp\Big\{\frac t{\alpha(t)^2}\Hcal_{\sM}^{\ssup{t}}(\ell_t)+\Rcal_{\sM}^{\ssup {t}}( \ell_t)
    \Big\}
    \one\{\supp(\ell_t)\subseteq B_{R\alpha(t)}\} \Big],$$
We will see that $\Hcal_{\sM}^{\ssup{t}}$ gives the main term
and $\Rcal_{\sM}^{\ssup {t}}$ a small remainder in
the limit $t\to\infty$, followed by $M\to\infty$. To separate
the two factors coming from this split, we use H\"older's
inequality. For any small $\eta>0$, we have
    \begin{equation}\label{hoelder}
    \begin{aligned}
    \big\langle u^{\xi_{t}}_{R\alpha(t)}(t,0) \big\rangle
    &\leq \me_0\Big[ \exp\Big\{ {(1+\eta)} \frac t{\alpha(t)^2}\Hcal_{\sM}^{\ssup{t}}(\ell_t) \Big\}
    \one\{\supp(\ell_t)\subseteq B_{R\alpha(t)}\} \Big]^\frac{1}{1+\eta}\\
    &\quad \times \me_0\Big[ \exp\Big\{ \sfrac{1+\eta}{\eta}\Rcal_{\sM}^{\ssup {t}}( \ell_t)\Big\}
    \one\{\supp(\ell_t)\subseteq B_{R\alpha(t)}\} \Big]^{\sfrac{\eta}{1+\eta}}.
    \end{aligned}
    \end{equation}
We show later that the second factor is asymptotically negligible, more precisely, we show that
\begin{equation}\label{remainderesti}
\limsup_{M\to\infty}\limsup_{t\to\infty}\frac {\alpha(t)^2}t\log\me_0\Big[
\exp\Big\{C\Rcal_{\sM}^{\ssup {t}}( \ell_t)\Big\}
\one\{\supp(\ell_t)\subseteq B_{R\alpha(t)}\} \Big]=0,\qquad \mbox{for }C>0.
\end{equation}
Let us first focus on the first term. Recall the definition of $\alpha(t)$ in \eqref{alphadef}
and the uniform convergence claimed in Proposition~\ref{variation}. For every $\eps>0$ and all
sufficiently large times $t$, we obtain the upper bound
\begin{align}
 \Hcal_{\sM}^{\ssup{t}}(\ell_t)&\leq \sfrac {\alpha(t)^2}t \kappa(\sfrac{t}{\alpha(t)^d})\, \rho
\sum_{x\in B_{R\alpha(t)}}  \sfrac{\ell_t(x)}{t/\alpha(t)^d} \log\big( \sfrac{\ell_t(x)}{t/\alpha(t)^d} \big)
\one\{ \ell_t(x) \leq  \sfrac{Mt}{\alpha(t)^d}\}
+ \eps \,  (2R)^d\alpha(t)^{d}\, \sfrac {\alpha(t)^2}t  \kappa(\sfrac{t}{\alpha(t)^d})\nonumber\\
&\leq \rho \sum_{x\in B_{R\alpha(t)}} \sfrac 1t \ell_t(x)\log\big( \sfrac{1}{t}\ell_t(x)\,\alpha(t)^d \big)
+ \eps \, (2R)^d = G_t(\sfrac 1t\ell_t)+ \eps \,  (2R)^d,\label{HMesti}
\end{align}
where we dropped the indicator, which we can do for $M\geq 1$ since $y\log y\geq 0$ for $y> M$,
and let\begin{equation}\label{Gdef}
G_t(\mu)=\rho\,\sum_{x\in B_{R\alpha(t)}}\mu(x)\log\big(\alpha(t)^d\mu(x)\big),
\qquad \mbox{ for }\mu\in\Mcal(\Z^d).
\end{equation}

The further analysis makes crucial use of an inequality derived in~\cite{BHK04}.
In \cite{BHK04}, the law of the local times are investigated,
and an explicit formula is derived for the density of the
local times on the range of the random walk. This
explicit formula makes it possible to give strong upper bounds
on exponential functionals:

\begin{prop}
\label{prop-BHK04}
For any
finite set $B\subseteq \Z^d$ and any measurable
functional $F\colon \Mcal_1(B)\to \R$,
    \begin{equation}\label{supertool}
    \me_0\Big[ e^{tF(\frac 1t\ell_t)} \, \one\{\supp(\ell_t)\subseteq B\} \Big]
    \leq \exp\Big\{ t\sup_{\mu\in\Mcal_1(B)}\! \Big[
    F(\mu) - \sfrac{1}{2} \sum_{x \sim y} \big( \sqrt{\mu(x)}-\sqrt{\mu(y)} \big)^2
    \Big] \Big\}\,(2dt)^{|B|}|B|.
    \end{equation}
\end{prop}

We substitute \eqref{HMesti} into \eqref{hoelder} and apply \eqref{supertool} for $F=(1+\eta)G_t/\alpha(t)^2$
and $B=B_{R\alpha(t)}$ and note that $(2dt)^{|B_{R\alpha(t)}|}\leq e^{o(t/\alpha(t)^2)}$.
Hence, we obtain that the first term on the right hand side of
\eqref{hoelder} can be estimated by
    \begin{equation}\label{toolappl}
    \begin{aligned}
    \me_0\Big[ &\exp\Big\{ {(1+\eta)} \frac t{\alpha(t)^2}\Hcal_{\sM}^{\ssup{t}}(\ell_t) \Big\}
    \one\{\supp(\ell_t)\subseteq B_{R\alpha(t)}\} \Big]\\
    &\leq e^{o(t/\alpha(t)^2)}\exp\Big\{-t\Big[ \chi^{\rm d}\big(\sfrac{\widetilde\rho}{\alpha(t)^2}\big)-
    \sfrac {1}{\alpha(t)^2}\, \big( d\widetilde\rho\log\alpha(t) + \eps\,(2R)^d \big) \Big]\Big\},
    \end{aligned}
    \end{equation}
where we abbreviated $\widetilde\rho=(1+\eta)\rho$ and introduced
    \begin{equation}\label{varGH}
    \chi^{\rm d}(\delta)=\inf_{\mu\in\Mcal_1(\Z^d)}\Big[\sfrac{1}{2} \sum_{x \sim y} \big( \sqrt{\mu(x)}-
    \sqrt{\mu(y)} \big)^2-\delta \sum_{x\in\Z^d}\mu(x)\log\mu(x)\Big],\qquad \mbox{ for }\delta>0,
    \end{equation}
the discrete variant of $\chi(\rho)$ in \eqref{chidef}, which
was studied in G\"artner and den Hollander~\cite{GH99}. In Proposition~3 and
the subsequent remark they show that
    $$
    \chi^{\rm d}(\delta)= \frac{d\delta}{2}\, \Big( \log\frac{\pi e^2}{\delta} + o(1) \Big), \qquad
    \mbox{ as } \delta\downarrow 0.
    $$
Substituting this into \eqref{toolappl}, we obtain
    \begin{equation}\label{toolappl1}
    \begin{aligned}
    \limsup_{t\to\infty}&\frac {\alpha(t)^2}t\log\me_0\Big[ \exp\Big\{ {(1+\eta)}
    \frac t{\alpha(t)^2}\Hcal_{\sM}^{\ssup{t}}(\ell_t) \Big\}
    \one\{\supp(\ell_t)\subseteq B_{R\alpha(t)}\} \Big]\\
    &\leq-\frac{\widetilde \rho d}2\log\frac{\pi e^2}{\widetilde\rho} + \eps (2R)^d
    =-\chi(\widetilde \rho)+\eps (2R)^d,
    \end{aligned}
    \end{equation}
as can be seen from Proposition~\ref{lem-chiident}. Using \eqref{toolappl1} together
with \eqref{remainderesti} in (\ref{hoelder}) and letting $M\to\infty$, $\eps\downarrow0 $
and $\eta\downarrow0 $,  gives the desired upper bound and finishes the proof of
Proposition~\ref{largedeviat1}(i) subject to the proof
of~\eqref{remainderesti}.

It remains to investigate the second term in (\ref{hoelder}),
i.e., to prove \eqref{remainderesti}. We first estimate
$\Rcal_{\sM}^{\ssup t}(\ell_t)$ (recall \eqref{remainder})
from above in terms of a nice functional of $\ell_t$. Since
we have to work uniformly for arbitrarily large local times, it
is not possible to estimate against a functional of the form
$\sum_x\ell_t(x)\log \ell_t(x)$, but we succeed
in finding an upper bound of the form $(\sum_x\ell_t(x)^{q})^{1/q}$
for some $q>1$ close to 1. Then Proposition~\ref{cutting} can be applied
and yields \eqref{remainderesti}.

We fix $\delta\in(0,\frac 12]$ and note that there exist $A>1$, $t_0>0$ such that
    \begin{equation}
    \label{Hinftyesti}
    \frac{H(ty)-yH(t)}{\kappa(t)} \leq A y^{1+\delta^2/3}
    \qquad\mbox{ for any $y\geq 1$ and $t>t_0$.}
    \end{equation}
Indeed, this follows from \cite[Theorem 3.8.6(a)]{BGT87}. Therefore,
we obtain that
    \eq\label{Hinftyesti1}
    H(\ell_t(x))-\ell_t(x) \sfrac{\alpha(t)^d}{t} H( \sfrac{t}{\alpha(t)^d})
    \leq A\kappa(t\alpha(t)^{-d}) \Big(\frac{\ell_t(x)}{t\alpha(t)^{-d}}\Big)^{1+\delta^2/3}.
    \en
We pick now $\eps>0$ such that
    \begin{equation}
    \label{epsdef}
    1+\delta^2/3-\eps=1/(1+\delta).
    \end{equation}
For any $\mu\in\Mcal_1(\Z^d)$, we use Jensen's inequality together with (\ref{epsdef})
as follows:
    \begin{equation}
    \label{estiHinfty2}
    \begin{aligned}
    \sum_{x\colon \mu(x)>M}
    \mu(x)^{1+\delta^2/3}&=\Bigl(\sum_{x\colon \mu(x)>M}
    \mu(x)^\eps\Bigr)\,\sum_{x\colon \mu(x)>M}
    \frac{\mu(x)^\eps}{\sum_{x\colon \mu(x)>M}
    \mu(x)^\eps} \mu(x)^{1+\delta^2/3-\eps}\\
    &\leq \Bigl(\sum_{x\colon \mu(x)>M}
    \mu(x)^\eps\Bigr)\,\Bigl(\sum_{x\colon \mu(x)>M}
    \frac{\mu(x)^{1+\eps}}{\sum_{x\colon \mu(x)>M}
    \mu(x)^\eps}\Bigr)^{\frac 1{1+\delta}}\\
    &=\Bigl(\sum_{x\colon \mu(x)>M}
    \mu(x)^\eps\Bigr)^{1-\frac 1{1+\delta}}\Bigl(\sum_{x\colon \mu(x)>M}
    \mu(x)^{1+\eps}\Bigl(\frac {\mu(x)}M\Bigr)^{\delta-\eps}\Bigr)^{\frac 1{1+\delta}}\\
    &\leq M^{\frac{\delta}{1+\delta}(\eps-1)}
    M^{\frac{\eps-\delta}{1+\delta}}\Bigl(\sum_{x}
    \mu(x)^{1+\delta}\Bigr)^{\frac 1{1+\delta}}
    =M^{\eps-\frac{2\delta}{1+\delta}}\Bigl(\sum_{x}
    \mu(x)^{1+\delta}\Bigr)^{\frac 1{1+\delta}},
    \end{aligned}
    \end{equation}
where we used in the last step that in the first integral on the right, $\mu^\eps\leq M^{\eps-1} \mu$
on $\{\mu>M\}$, and hence the first term
on the right is not bigger than one, as the exponent is positive and $\mu\in\Mcal_1(\Z^d)$.
We write $q=1+\delta$. We apply the above to $\mu=\frac 1t\ell_t$ and $M$
replaced by $\sfrac{M}{\alpha(t)^d}$, to obtain that
    \begin{equation}\label{Hinftyest2}
    \begin{aligned}
    \sum_{x} \Big(\frac{\ell_t(x)}{t\alpha(t)^{-d}}\Big)^{1+\delta^2/3}
    \one\{ \ell_t(x) >  \sfrac{Mt}{\alpha(t)^d} \}
    &\leq \alpha(t)^{d(1+\delta^2/3)}\big( \sfrac{M}{\alpha(t)^d}\big)^{\eps-\frac{2\delta}{1+\delta}}
    t^{-1}\|\ell_t\|_q\\
    &= M^{\eps-\frac{2\delta}{1+\delta}} \alpha(t)^{d(1+\frac{\delta}{1+\delta})}
    t^{-1}\|\ell_t\|_q.
    \end{aligned}
    \end{equation}
We recall \eqref{remainder}, use \eqref{Hinftyesti1} and the definition of $\alpha(t)$
in (\ref{alphadef}). With the help of \eqref{Hinftyest2} we arrive at
$$  \begin{aligned}
    \Rcal_{\sM}^{\ssup t}(\ell_t)&\leq \kappa(\sfrac{t}{\alpha(t)^d})
    \sum_{x} \Big(\frac{\ell_t(x)} {t\alpha(t)^{-d}}\Big)^{1+\delta^2/3}
    \one\{ \ell_t(x) > A \sfrac{Mt}{\alpha(t)^d} \}\\
    &\leq AM^{\eps-\frac{2\delta}{1+\delta}} \alpha(t)^{-(2+d)+d(1+\frac{\delta}{1+\delta})} \|\ell_t\|_q\\
    &= AM^{\eps-\frac{2\delta}{1+\delta}} \alpha(t)^{-\frac 1q [d+(2-d)q]} \|\ell_t\|_q,
    \end{aligned}$$
where we recall that $q=1+\delta$ and therefore $-(2+d)+d(1+\frac{\delta}{1+\delta})=-\frac 1q [d+(2-d)q]$.
Put $\theta=AM^{\eps-\frac{2\delta}{1+\delta}}$, and observe that $\theta\downarrow 0$ as $M\uparrow \infty$
for $\delta>0$ small enough, since
$\eps-\frac{2\delta}{1+\delta} = \frac{\delta^2}{3}-\frac{\delta}{1+\delta}<0$
for $\delta>0$ small enough.
Hence, \eqref{remainderesti} follows immediately from Proposition~\ref{cutting}.
This completes the proof of Proposition~\ref{largedeviat1}(i).

\subsection{Proof of Proposition~\ref{largedeviat1}(ii)}
\label{sec-pfprops2}

\noindent Recall from \eqref{shirescpot} the rescaled version, $\overline\xi_t$,
of the vertically shifted potential, $\xi_t$, defined in \eqref{vtdef}. Furthermore,
introduce the normalised, scaled version of the random walk local times,
    $$
    L_t(x) := \frac{\alpha(t)^d}{t} \, \ell_t\big( \lfloor x \alpha(t) \rfloor
    \big) , \qquad \mbox{ for } x \in \R^d,
    $$
and note that $L_t$ is an $L^1$-normalised random step function. Note that $\supp(L_t)\subseteq Q_{\sR}$
if $\supp(\ell_t)\subseteq B_{R\alpha(t)}$ where we abbreviated $Q_{\sR}=[-R,R]^d$. We start from \eqref{utident2}.
Let
    $$
    \widehat{H}_t(y)= \frac{H\big(y\frac{t}{\alpha(t)^d}\big)-
    y H\big(\frac{t}{\alpha(t)^d}\big)}{ \kappa\big(\frac{t}{\alpha(t)^d}\big)},\qquad \mbox{for }t,y>0,
    $$
and recall that $\widehat H_t$ converges to $\widehat H$, uniformly on all compact sets. Now
the exponent on the right hand side of \eqref{utident2} can be rewritten as follows.
    $$
    \begin{aligned}
    -\alpha(t)^d H\big(\sfrac{t}{\alpha(t)^d}\big) +
    \sum_{z\in B_{R\alpha(t)}} H(\ell_t(z))
    & = -\alpha(t)^d \int_{Q_{\sR}}L_t(x)\,H\big(\sfrac{t}{\alpha(t)^d}\big)\,dx +
    \alpha(t)^d \int_{Q_{\sR}} H\big( \sfrac{t}{\alpha(t)^d} L_t(x)\big) \, dx \\
    & =\alpha(t)^d\kappa\big({\textstyle{\frac{t}{\alpha(t)^d}}}\big) \int_{Q_{\sR}}\widehat{H}_t\big(L_t(x)\big) \, dx =\sfrac{t}{\alpha(t)^2} \skrih^{\ssup{t}}_{\sR}(L_t),
    \end{aligned}
    $$
where we use the definition of $\alpha(t)$ in \eqref{alphadef} and introduce the functional
    $$
    \skrih^{\ssup{t}}_{\sR}(f) = \int_{Q_{\sR}} \widehat H_t\big(f(x)\big)\, dx.
    $$
Hence,
    \begin{equation}\label{utident3}
    \big\langle u^{\xi_{t}}_{R\alpha(t)}(t, 0) \big\rangle
    = \me_0\bigg[\exp\Big( \frac{t}{\alpha(t)^2}\skrih^{\ssup t}_{\sR}(L_t)\Big)
    \one\{\supp(L_t)\subseteq Q_{\sR}\} \bigg].
    \end{equation}
A key ingredient in the proof of Proposition~\ref{largedeviat1}(ii) is the
large deviation principle for $(L_t \colon t>0)$ as formulated in the following
proposition:

\begin{prop}
\label{prop-GKS04}
Fix $R>0$. Under $\P_0\{\,\cdot\, ,\supp(L_t)\subseteq Q_{\sR}\}$,
the rescaled local times process $(L_t \colon t>0)$ satisfies a
large deviation principle as $t\uparrow\infty$ on the set of
$L^1$-normalized functions $Q_{\sR}\to\R$,
equipped with the weak topology induced by test integrals
against all continuous functions, where the speed of the large
deviation principle is $t\alpha(t)^{-2}$,
and the rate function is $g^2\mapsto\|\nabla g\|_2^2$, on the set of all
$g\in H^1(\R^d)$ with $\supp(g)\subseteq Q_{\sR}$, and is equal to
$\infty$ outside this set.
\end{prop}

\begin{Proof}{Proof.}
This large deviation principle is stated in \cite[Lemma~3.2]{GKS04}
in the discrete-time case, and is proved in \cite[Section 6]{GKS04}.
The proof in the continuous-time case is very similar.
\qed
\end{Proof}

In order to apply the large deviation principle in Proposition \ref{prop-GKS04}
to obtain a lower bound for the right hand side of
 \eqref{utident3}, we need the lower-bound half of Varadhan's lemma, and we have to replace
$\skrih^{\ssup{t}}_{\sR}$ by its limiting version
    \begin{equation}\label{HRdef}\skrih_{\sR}(f)
    =\rho\int_{Q_{\sR}}f(x)\log f(x)\, dx.
    \end{equation}
However, the latter is technically not so easy. Inserting the indicator on the event $\{\|L_t\|_\infty< M\}$
for any $M>1$ would make it possible to use the locally uniform convergence of $\widehat H_t(y)$ towards
$\rho y\log y$, but this event is not open in the topology of the large deviation principle.
Therefore, similarly to the proof of the upper bound, we have to split
$\skrih^{\ssup{t}}_{\sR}(L_t)$ into the sum of $\skrih_{\sR}(L_t)$ and a remainder term,
separate these two from each other by the use of H\"older's
inequality and apply Proposition~\ref{cutting} to the remainder term.
Let us turn to the details.

Since $H$ is convex with $H(0)=0$, we have $H(yt)\geq y H(t)$ for all $t>0$ and all $y\geq 1$. Therefore, $\widehat H_t(f(x))\geq 0$ on $\{x\colon f(x)>M\}$ for any $M>1$. Hence, we may estimate
    $$
    \begin{aligned}
    \skrih^{\ssup{t}}_{\sR}(f)&\geq \int_{Q_{\sR}} \1\{f(x)\leq M \}\widehat H_t\big(f(x)\big)\, dx=
    \rho \int_{Q_{\sR}} \1\{f(x)\leq M \}f(x)\log f(x)\, dx+o(1)\\
    &=\skrih_{\sR}(f)-\rho\int_{Q_{\sR}} \1\{f(x)> M \}f(x)\log f(x)\, dx+o(1).
    \end{aligned}
    $$
The remainder can be estimated, for any $\delta>0$, as follows. For any $f\colon Q_{\sR}\to[0,\infty)$ satisfying $\int f=1$,
    $$
    \begin{aligned}
    \int_{f>M}f\log f&=\frac 2\delta \Big(\int_{f>M}f\Big)\int_{f>M}\frac f{\int_{f>M}f}
    \log f^{\delta/2}\leq \frac 2\delta \Big(\int_{f>M}f\Big)\log\frac{\int_{f>M}f^{1+\delta/2}}{\int_{f>M}f}\\
    &\leq \frac 2\delta \Big(\int_{f>M}f\Big)\log\frac{M^{-\delta/2}\int_{f>M}f^{1+\delta}}{\int_{f>M}f}\leq
    \frac 2\delta \Big(\int_{f>M}f\Big)\Big(\frac{M^{-\delta/2}
    \int_{f>M}f^{1+\delta}}{\int_{f>M}f}\Big)^{\frac 1{1+\delta}}\\
    &=\frac 2\delta M^{-\frac \delta{2+2\delta}}\Big(\int_{f>M}f\Big)^{\frac \delta{1+\delta}}
    \|f\|_q\leq \frac 2\delta M^{-\frac \delta{2+2\delta}}\|f\|_q,
    \end{aligned}
    $$
where we put $q=1+\delta$.
Altogether, we have, abbreviating $\theta=2\frac\rho\delta M^{-\frac \delta{2+2\delta}}$,
    \begin{equation}\label{utident4}
    \big\langle u^{\xi_{t}}_{R\alpha(t)}(t, 0) \big\rangle
    \geq \me_0\bigg[\exp\Big( \frac{t}{\alpha(t)^2}\big(\skrih_{\sR}(L_t)-\theta \|L_t\|_q\big)\Big)
    \one\{\supp(L_t)\subseteq Q_{\sR}\} \bigg]e^{o(t/\alpha(t)^2)}.
    \end{equation}
Similarly to the proof of the upper bound, the main contribution will turn out to come from $\skrih_{\sR}$,
and the $q$-norm is a small remainder. In order to separate the two from each other, we use H\"older's
inequality to estimate, for some small $\eta>0$,
    \begin{equation}\label{Hoelderprel}
    \begin{aligned}
    &\me_0\bigg[\exp\Big( \frac{t}{\alpha(t)^2}(1-\eta)\skrih_{\sR}(L_t)\Big)
    \one\{\supp(L_t)\subseteq Q_{\sR}\} \bigg]\\
    &\quad \leq \me_0\bigg[\exp\Big( \frac{t}{\alpha(t)^2}\big(\skrih_{\sR}(L_t)-\theta \|L_t\|_q\big)\Big)
    \one\{\supp(L_t)\subseteq Q_{\sR}\} \bigg]^{1-\eta}\\
    &\qquad \times
    \me_0\bigg[\exp\Big( \frac{t}{\alpha(t)^2}\frac {1-\eta}\eta\theta \|L_t\|_q\Big)
    \one\{\supp(L_t)\subseteq Q_{\sR}\} \bigg]^{\eta}.
    \end{aligned}
    \end{equation}
This effectively yields a lower bound on the expected value in \eqref{utident4}
of the form
    \begin{equation}\label{Hoelder}
    \begin{aligned}
    \big\langle u^{\xi_{t}}_{R\alpha(t)}(t, 0) \big\rangle
    &\geq\me_0\bigg[\exp\Big( \frac{t}{\alpha(t)^2}(1-\eta)\skrih_{\sR}(L_t)\Big)
    \one\{\supp(L_t)\subseteq Q_{\sR}\} \bigg]^{\frac 1{1-\eta}}\\
    &\quad\times\me_0\bigg[\exp\Big( \frac{t}{\alpha(t)^2}\frac {1-\eta}\eta\theta \|L_t\|_q\Big)
    \one\{\supp(L_t)\subseteq Q_{\sR}\} \bigg]^{-\frac \eta{1-\eta}}e^{o(t/\alpha(t)^2)}.
    \end{aligned}
    \end{equation}
From Proposition~\ref{cutting} it follows that the second expectation on the right is negligible
in the limit $t\to\infty$, followed by $M\to\infty$, i.e., $\theta\downarrow 0$. Hence, we
can concentrate on the first term.  To apply the lower-bound half of Varadhan's lemma,
see \cite[Lemma~4.3.4]{DZ98}, we need the following lower semi-continuity property of the
function $\skrih_{\sR}$:

\begin{lemma}\label{lem-semic}
Let $f\colon Q_{\sR}\to[0,\infty)$ be continuous. Then  $\skrih_{\sR}$ is lower semi-continuous in $f$
in the topology induced by pairing with all continuous functions $Q_{\sR}\to[0,\infty)$.
\end{lemma}

\begin{Proof}{Proof.} Let $(f_n \colon n\in\N)$ be a family in $L^1(Q_{\sR})$ such that
$\langle f_n,\psi\rangle\to \langle f,\psi\rangle$ as $n\to\infty$ for any continuous function
$\psi\colon Q_{\sR}\to \R$. We have to show that $\liminf_{n\to\infty}\skrih_{\sR}(f_n)\geq \skrih_{\sR}(f)$.

For any $s\in(0,\infty)$ we denote by $g_s$ the tangent to $y \mapsto \phi(y):=\rho y \log y$ in $s$,
i.e., $g_s(y) = \rho(1+\log s) y - \rho s, \mbox{ for all } y\in \R.$
By convexity we have $g_s\leq \phi$ for any $s\in(0,\infty)$. Therefore, for any $0<\eps<1/e$,
    $$
    \skrih_{\sR}(f_n) = \int_{Q_{\sR}}\phi\big(f_n(x)\big) \, dx\geq\int_{Q_{\sR}} g_{f(x)\lor \eps}\big(f_n(x)\big) \, dx
    =\rho \big\langle 1+\log (f\lor \eps), f_n \big\rangle
    - \rho  \big\langle f \lor \eps, f_n\big\rangle.
    $$
Letting $n\to\infty$, we obtain, using the boundedness and continuity of $\log (f\lor \eps)$,
    $$
    \begin{aligned}
    \liminf_{n\to\infty}\skrih_{\sR}(f_n)&\geq  \rho \big\langle 1+\log (f\lor \eps), f \big\rangle
    - \rho  \big\langle f \lor \eps, f \big\rangle\\
    & \geq \rho \, \int_{Q_{\sR}} f(x) \log f(x) \one_{\{ f(x) > \eps \}} \, dx + \int_{Q_{\sR}}  g_\eps(f(x)) \one_{\{ f(x) \leq \eps \}} \, dx \\
    & \geq
    \rho \, \int_{Q_{\sR}} f(x) \log f(x)
    \, dx + \int_{Q_{\sR}} \big( f(x)(1+\log \eps)-\eps\big) \one_{\{ f(x) \leq \eps \}} \, dx.
    \end{aligned}
    $$
The second summand is bounded from below by Leb$(Q_{\sR}) \eps \log \eps$, which converges
to zero as $\eps\downarrow 0$. This completes the proof.
\qed\end{Proof}

Now we can apply \cite[Lemma~4.3.4]{DZ98} and obtain
    $$
    \begin{aligned}
    \liminf_{t\to\infty}&\frac{\alpha(t)^2}t\log \me_0\bigg[\exp\Big( \frac{t}{\alpha(t)^2}(1-\eta)
    \skrih_{\sR}(L_t)\Big)\one\{\supp(L_t)\subseteq Q_{\sR}\} \bigg]\\
    &\geq -\inf\Big\{\|\nabla g\|_2^2-(1-\eta)\skrih_{\sR}(g^2)
    \colon g\in H^1(\R^d)\cap\Ccal(\R^d),\|g\|_2=1,\supp(g)\subseteq Q_{\sR}\Big\}.
    \end{aligned}
    $$
Letting $\eta\downarrow0$ and $R\uparrow\infty$, it is easy to see that the right hand side tends to
$-\chi(\rho)$ defined in \eqref{chidef}. Indeed, use appropriate continuous cut-off versions
$g_{(1-\eta)\rho}^{\ssup R}$ of the minimiser $g_{(1-\eta)\rho}$ in \eqref{Gauss} to verify this claim.
Using this on the right hand side of \eqref{Hoelder} and recalling Proposition~\ref{cutting}, we see
that the proof of the lower bound in Proposition~\ref{largedeviat1}(ii) is finished.

\section{The almost-sure asymptotics: Proof of Theorem~\ref{asasy}}\label{sec-proofquen}

\noindent We again derive upper and lower bounds, following the strategy in
\cite[Section 5]{BK01}. Recall the scale function $\beta(t)$ defined in \eqref{betatdef} and let
    \begin{equation}\label{xibetadef}
    \xi_{\beta(t)}(z)=\xi(z)-H\Big(\frac{\beta(t)}{\alpha(\beta(t))^{d}}\Big)\frac{\alpha(\beta(t))^{d}}
    {\beta(t)}
    \end{equation}
denote the appropriately vertically shifted potential (compare to \eqref{vtdef}). Then Theorem~\ref{asasy}
is equivalent to the assertion
\begin{equation}\label{mainaimquen}
\lim_{t\uparrow\infty}\frac{\alpha(\beta(t))^2}{t}\log u^{\xi_{\beta(t)}}(t,0)=-\widetilde\chi(\rho),\qquad
\mbox{almost surely,}
\end{equation}
where $\widetilde\chi(\rho)=\rho(d-\sfrac d2 \sfrac{\rho}{\pi } + \log \sfrac\rho{e})
=-\sup\{\lambda(\psi)\colon\psi\in\Ccal(\R^d),\Lcal(\psi)\leq 1\}$, see Section~\ref{asvarform}.

\subsection{Proof of the upper bound in {\bf{(\ref{mainaimquen})}}}

\noindent Let $r(t)=t \log t$ and apply Lemma~\ref{compact} with $V=\xi_{\beta(t)}$ and with $R$
replaced by $R\alpha(\beta(t))$. Furthermore, take logarithms, multiply with $\alpha(\beta(t))^2/t$
and let $t\uparrow\infty$. As in \eqref{firsttermneg}, one shows that the first term is negligible.
Hence, we obtain that
$$    \limsup_{t\uparrow\infty}\frac {\alpha(\beta(t))^2}{t}\log u^{\xi_{\beta(t)}}(t,0)\leq
    \frac C{R^2}+\limsup_{t\uparrow\infty}\Bigl[\alpha(\beta(t))^2\max_{z\in B(t)}
    \lambda_{z;2R\alpha(\beta(t))}(\xi_{\beta(t)})\Big],$$
where $B(t)=B_{R\alpha(\beta(t))}(t)$ (recall the definition
$B_{\sR}(t)=B_{r(t)+\lfloor 2R\rfloor}$ from Lemma~\ref{compact}).
Let $(\lambda_i(t) \colon i=1,\dots,N(t))$, with $N(t)=|B_{\sR}(t)|$,
be a deterministic enumeration of the random variables
$\lambda_{z;2R\alpha(\beta(t))}(\xi_{\beta(t)})$ with
$z\in B(t)$. Note that these random variables are identically
distributed (but not independent) and that, by \eqref{fromBK} and
Proposition~\ref{largedeviat1}(i), their exponential moments are estimated by
    \begin{equation}\label{expmomtildexi}
    \limsup_{t\uparrow\infty}\frac {\alpha(\beta(t))^2}{\beta(t)}\log\big\langle
    e^{\beta(t)\, \lambda_1(t)}\big\rangle\leq -\chi(\rho).
    \end{equation}
We next show that, for any $\eps>0$, almost surely,
    \begin{equation}\label{aimquen1}
    \limsup_{t\uparrow\infty}\alpha(\beta(t))^2\max_{i=1}^{N(t)}\lambda_i(t)
\leq -{\widetilde \chi}(\rho)+\eps,
    \end{equation}
which completes the proof of the upper bound in \eqref{mainaimquen}.

To prove \eqref{aimquen1}, one first realizes that is suffices to show \eqref{aimquen1}
only for $t\in\{e^n\colon n\in\N\}$, since the functions $t\mapsto \alpha(t)$,
$t\mapsto \beta(t)$, and $t \mapsto  H(t)/t$ are slowly varying, and $t \mapsto N(t)$,
$R\mapsto \lambda_{R}(\xi_{\beta(t)})$ are increasing. Let
$$p_n=\Prob\Big\{\max_{i=1}^{N(e^n)} \lambda_i(e^n)\geq \frac{-{\widetilde \chi(\rho)}+\eps}
{\alpha(\beta(e^n))^2}\Big\}.$$
We recall that $\beta(e^n)\alpha(\beta(e^n))^{-2} \sim dn$. Using Chebyshev's inequality
and~\eqref{expmomtildexi}, we estimate, for any $k>0$,
    \begin{equation}\label{chebyappl}
    \begin{aligned}
    p_n&\leq N(e^n)\Prob\Big\{e^{k\beta(e^n)\lambda_1(e^n)}\geq e^{-k\beta(e^n)\alpha(\beta(e^n))^{-2}
    ({\widetilde \chi}(\rho)-\eps)}\Big\}\\
    &\leq e^{n(d+o(1))} e^{nkd(\widetilde\chi(\rho)-\eps)}
    \Big\langle e^{k\beta(e^n)\lambda_1(e^n)}\Big\rangle.
    \end{aligned}
    \end{equation}
In order to evaluate the last expectation, we intend to apply \eqref{expmomtildexi} with
$\beta(t)$ replaced by $k\beta(t)$. For this purpose, we note that we can replace
$\alpha(\beta(t))$ by $\alpha(k \beta(t))$ in \eqref{expmomtildexi}, since $\alpha$
is slowly varying. Also,
    $$\begin{aligned}
    k \beta(t) \lambda_{R\alpha(k \beta(t))}(\xi_{\beta(t)})&  =
    k \beta(t) \lambda_{R\alpha(k \beta(t))}(\xi_{k\beta(t)})
    - k \beta(t) \Big[ \frac{H(\beta(t) \alpha(\beta(t))^{-d})}{\beta(t) \alpha(\beta(t))^{-d}}
    -  \frac{H(k\beta(t) \alpha(k\beta(t))^{-d})}{k\beta(t) \alpha(k\beta(t))^{-d}} \Big],
    \end{aligned}$$
where we use that by \eqref{xibetadef}, the field $\xi_{\beta(t)}-\xi_{k\beta(t)}$
is constant and deterministic. Now we use \eqref{basic} and \eqref{alphadef}, to see that
the deterministic term is equal to
$$\begin{aligned}
 k \beta(t) \Big[ & \frac{H(\beta(t) \alpha(\beta(t))^{-d})}{\beta(t) \alpha(\beta(t))^{-d}}
-  \frac{H(k\beta(t) \alpha(k\beta(t))^{-d})}{k\beta(t) \alpha(k\beta(t))^{-d}} \Big] \\
& = \alpha(\beta(t))^d \, \Big( k \,  H(\beta(t) \alpha(\beta(t))^{-d})
-H(k\beta(t) \alpha(\beta(t))^{-d}) \Big) + o(n)\\
& = -\alpha(\beta(t))^d \, \big(\widehat{H}(k)+o(1)\big) \, 
\kappa\big(\beta(t) \alpha(\beta(t))^{-d}\big) \,  + o(n)\\
& =  - \frac{\beta(t)}{\alpha(\beta(t))^2}\, \big(\rho k \log k + o(1)\big) (1 + o(1))\,  + o(n) \\
& =  - nd \, (\rho k \log k) (1 + o(1)).  
\end{aligned}$$
Hence,
    $$\Big\langle e^{k\beta(e^n)\lambda_1(e^n)}\Big\rangle
    \leq \exp\Big\{ -nd \big( k \chi(\rho) -  \rho k \log k + o(1)\big) \Big\}.
    $$
Using this in \eqref{chebyappl}, we arrive at
$$p_n \leq  \exp\Big\{ nd \big( 1 + k(\widetilde\chi(\rho)-\eps) -
k \chi(\rho) +  \rho k \log k + o(1)\big) \Big\}.$$
Choosing $k=\frac 1{\rho}$, we see that $p_n\leq e^{-nd(k\eps + o(1))}$. This
is summable over $n\in\N$, and the Borel-Cantelli lemma yields
that \eqref{aimquen1} holds almost surely. This completes the proof of the upper bound
in (\ref{mainaimquen}).

\subsection{Proof of the lower bound in {\bf{(\ref{mainaimquen})}}}

\noindent Our proof of the lower bound in (\ref{mainaimquen})
follows the strategy of \cite[Sect.~5.2]{BK01}. First we establish
that, with probability one, for any sufficiently large $t$, there is,
inside a `macrobox' of radius roughly~$t$, centred at the origin,
some `microbox' of radius $R\alpha(\beta(t))$ in which the random
field $\xi_{\beta(t)}$ has some shape with optimal spectral
properties. Then we obtain a lower bound for the Feynman-Kac
formula in \eqref{FK} by requiring that the random walk moves
quickly to that box and stays there for approximately~$t$ time
units. As a result, the contribution from that strategy is basically
given by the largest eigenvalue of $\Delta^{\rm d}+\xi$ in
that microbox. Rescaling and letting $R\uparrow\infty$, the lower
bound is derived from this.

Let us go to the details. We pick an increasing auxiliary scale function
$t\mapsto\gamma_t$ satisfying
    \begin{equation}\label{gammachoice}
    \begin{aligned}
    & \gamma_t=t^{1-o(1)},\qquad t-\gamma_t=t(1+o(1)),\\\
    &\gamma_t =o\Big(\frac t{\alpha(\beta(t))^{2}}\Big), \qquad
    \gamma_t \frac {H(\beta(t)\alpha(\beta(t))^{-d})}{\beta(t)\alpha(\beta(t))^{-d}}=
    o\Big(\frac t{\alpha(\beta(t))^{2}}\Big).
    \end{aligned}
    \end{equation}
(Note that the second requirement follows from the third.)
For example, $\gamma_t=t\alpha(\beta(t))^{-2}\eps_t$ with some
suitable $\eps_t\downarrow 0$ as a small inverse power of
$\log{t}$ satisfies \eqref{gammachoice}. This is obvious in the case where
$\lim_{s\uparrow\infty}H(s)/s=0$, and in the case where
$\lim_{s\uparrow\infty}H(s)/s=\infty$, it is also clear
since $H(s)/s$ diverges only subpolynomially in $s$, while
$\beta(t)=(\log{t})^{1+o(1)}$ and $\alpha$ is slowly varying.

The crucial step is to show that, in the \lq macrobox\rq\ $B_{\gamma_t}$, we find an appropriate
\lq microbox\rq. To fix some notation, let $Q_{\sR}=[-R,R]^d$ and let ${\mathcal C}(Q_{\sR})$ denote
the set of continuous functions $Q_{\sR}\to\R$. We need finite-space versions of the functionals
$\mathcal H, \mathcal L$ and $\lambda$ defined in \eqref{calHdef} and \eqref{Legendre}.
Recall the definition of ${\mathcal H}_{\sR}$ from \eqref{HRdef} and define its
Legendre transform ${\mathcal L}_{\sR}\colon \skric(Q_{\sR})\to(-\infty,\infty]$ by
    \begin{equation}\label{calLRdef}
    {\mathcal L}_{\sR}(\psi)=\sup\bigl\{\langle f,\psi \rangle-{\mathcal
    H}_{\sR}(f)\colon f\in \skric(Q_{\sR}),\, f\geq 0,\,\supp f\subseteq\supp\psi\bigr\}.
    \end{equation}
As in the proof of Proposition~\ref{lem-chiident} one can see that
$f=e^{\psi/\rho-1}$ is the unique maximizer in \eqref{calLRdef} with
$${\mathcal L}_{\sR}(\psi)= \frac\rho e \int _{Q_{\sR}} e^{\psi(x)/\rho }\, dx.$$

\begin{prop}[Existence of an optimal microbox]\label{microopt}
Fix $R>0$ and let $\psi\in\Ccal(Q_{\sR})$ satisfy $\Lcal_{\sR}(\psi)< 1$. Let $\eps>0$. Then, with
probability one, there exists $t_0>0$, depending also on $\xi$, such that, for all $t>t_0$,
there is $y_t\in B_{\gamma_t}$, depending on $\xi$, such that
    \begin{equation}\label{xilowbound}
    \xi_{\beta(t)}(y_t+z)\geq \frac 1{\alpha(\beta(t))^2}\psi
    \Big(\frac z{\alpha(\beta(t))}\Big)-\frac \eps{\alpha(\beta(t))^2},
    \qquad \mbox{ for }z\in B_{R\alpha(\beta(t))}.
    \end{equation}
\end{prop}
The proof of Proposition~\ref{microopt} is deferred to the end of this section.

Now we finish the proof of the lower bound in (\ref{mainaimquen}) subject to Proposition~\ref{microopt}.
Let $R,\eps>0$, and let $\psi\in\Ccal(Q_{\sR})$ be twice continuously
differentiable with $\Lcal_{\sR}(\psi)< 1$. Fix $\xi$ not belonging to the exceptional set of
Proposition~\ref{microopt}, i.e., let $t_0$ and $(y_t \colon t>t_0)$ in $B_{\gamma_t}$ be chosen such
that \eqref{xilowbound} holds for every~$t>t_0$. Fix $t>t_0$. In the Feynman-Kac formula
    $$
    u^{\xi_{\beta(t)}}(t,0)=\E_0\exp\Big\{\int_0^t \xi_{\beta(t)}(X(s))\,ds\Big\},
    $$
we obtain a lower bound by requiring that the random walk is at $y_t$ at time $\gamma_t$ and
remains within the microbox
$$B_{y_t,t}=y_t+B_{R\alpha(\beta(t))}$$
during the time interval $[\gamma_t, t]$. Using the Markov property at time $\gamma_t$, we obtain by this the lower bound
    \begin{equation}\label{uasylow1}
    \begin{aligned}
    u^{\xi_{\beta(t)}}(t,0)\geq \E_0 \Big[\exp & \Big\{\int_0^{\gamma_t}\xi_{\beta(t)}(X(s))\, ds\Big\}\delta_{y_t}(X(\gamma_t))\Big]\, \\ & \times \E_{y_t}\Big[ \exp\Big\{ \int_0^{t-\gamma_t}\xi_{\beta(t)}(X(s))\, ds\Big\}\1\{\tau_{y_t,t}>t-\gamma_t\}\Big],
    \end{aligned}
    \end{equation}
where $\tau_{y_t,t}=\inf\{s>0\colon X(s)\notin B_{y_t,t}\}$ denotes the exit time from the microbox
$B_{y_t,t}$. In the first expectation on the right side of \eqref{uasylow1}, we estimate $\xi$ from
below by its minimum $K=\essinf \xi(0)>-\infty$,
and in the second expectation we use \eqref{xilowbound} and shift spatially by $y_t$ to obtain
    \begin{equation}\label{uasylow2}
    \begin{aligned}
    u^{\xi_{\beta(t)}}(t,0)&\geq \exp\Big\{\gamma_t \Big[K-\frac {H(\beta(t)\alpha(\beta(t))^{-d})}
    {\beta(t)\alpha(\beta(t))^{-d}}\Big]\Big\}\P_0\{X(\gamma_t)=y_t\} \\
    &\quad \times
    e^{-\eps (t-\gamma_t)\alpha(\beta(t))^{-2}}
    \E_0\Big[ \exp\Big\{\int_0^{t-\gamma_t} \psi_t(X(s))\, ds \Big\}\1\{\tau_{0,t}>t-\gamma_t\} \Big] ,
    \end{aligned}
    \end{equation}
where we have denoted $\psi_t(\cdot)=\alpha(\beta(t))^{-2}\psi(\cdot\,\alpha(\beta(t))^{-1})$.
By our choice in \eqref{gammachoice}, the first term on the right side of \eqref{uasylow2} is
$e^{o(t\alpha(\beta(t))^{-2})}$. Now, by choosing a path from the origin to $y_t$ consisting
of $k$ steps for $k=\lfloor\gamma_t\rfloor$ or $k=\lfloor\gamma_t\rfloor+1$,
    $$
    \P_0\{X(\gamma_t)=y_t\}\geq \big( \sfrac{1}{2d} \big)^k \,
    \prob\{ \sigma(1) + \cdots + \sigma(k) \leq \gamma_t <  \sigma(1) + \cdots + \sigma(k+1) \},
    $$
where $\sigma(1), \sigma(2), \ldots$ are independent exponential random variables with
mean $1/2d$. Using that
    $$
    \prob\big\{ \sigma(1) + \cdots + \sigma(k) \leq \gamma_t <  \sigma(1) + \cdots + \sigma(k+1) \big\}
    \geq\prob\big\{ \sigma(1) + \cdots + \sigma(k) \in [\sfrac{\gamma_t}2, \gamma_t) \big\}
    \, \prob\big\{ \sigma(0) \geq \sfrac{\gamma_t}2 \big\},
    $$
and Cram\'er's theorem, we obtain the lower bound
    $$
    \P_0\{X(\gamma_t)=y_t\}\geq e^{-\Ocal(\gamma_t)}=e^{-o(t\alpha(\beta(t))^{-2})}.
    $$
By an eigenfunction expansion we have that
    $$
    \begin{aligned}
    \E_0\Big[ \exp\Big\{ & \int_0^{t-\gamma_t} \psi_t(X(s))\, ds \Big\}\1\{\tau_{0,t}>t-\gamma_t\}
    \Big] \\
    & \geq \E_0\Big[ \exp\Big\{ \int_0^{t-\gamma_t} \psi_t(X(s))\, ds \Big\}\1\{\tau_{0,t}>t-\gamma_t,
    X(t-\gamma_t)=0\} \Big] \\
    & \geq \exp\Big\{(t-\gamma_t)\lambda^{\rm d}(t)\Big\}{\rm e}_t(0)^2,
    \end{aligned}
    $$
where $\lambda^{\rm d}(t) $ is the principal eigenvalue of $\Delta^{\rm d}+\psi_t$ in the box
$B_{R\alpha(\beta(t))}$ with zero boundary condition, and ${\rm e}_t$ is the corresponding
positive $\ell^2$-normalized eigenvector. Summarising these estimates and recalling
from \eqref{gammachoice} that $t-\gamma_t=t(1+o(1))$, we obtain, almost surely,
    \begin{equation}\label{uasylow3}
    \liminf_{t\uparrow \infty}\frac{\alpha(\beta(t))^2}t\log u^{\xi_{\beta(t)}}(t,0)\geq -\eps+
    \liminf_{t\uparrow \infty}\alpha(\beta(t))^2 \lambda^{\rm d}(t)+\liminf_{t\uparrow \infty}
    \frac{\alpha(\beta(t))^2}t\log {\rm e}_t(0)^2.
    \end{equation}
We now define the continuous counterpart $\lambda_{\sR}$ of $\lambda^{\rm d}(t)$, which is the finite-space
version of the spectral radius defined in \eqref{Legendre}:
    \begin{equation}\label{lambdaRdef}
    \lambda_{\sR}(\psi) = \sup \big\{ \langle \psi, g^2 \rangle -\| \nabla g \|_2^2 \colon
    g\in H^1(\R^d), \,  \|g\|_2=1, \supp\, g \subseteq Q_{\sR} \big\}.
    \end{equation}
According to \cite[Lemma 5.3]{BK01},
    $$
    \liminf_{t\uparrow \infty}\alpha(\beta(t))^2 \lambda^{\rm d}(t)\geq
    \lambda_{\sR}(\psi)\qquad\mbox{and}\qquad
    \liminf_{t\uparrow \infty}\frac{\alpha(\beta(t))^2}t\log {\rm e}_t(0)^2\geq 0.
    $$
Using this in \eqref{uasylow3}, we obtain
    \begin{equation}\label{uasylow4}
    \liminf_{t\uparrow \infty}\frac{\alpha(\beta(t))^2}t\log u^{\xi_{\beta(t)}}(t,0)\geq -\eps+\lambda_{\sR}(\psi),
    \end{equation}
for any $\eps>0$ and for any twice continuously differentiable function $\psi\in\Ccal^2(Q_{\sR})$ satisfying
$\Lcal_{\sR}(\psi)<1$. Hence,
    $$
    \liminf_{t\uparrow \infty}\frac{\alpha(\beta(t))^2}t\log u^{\xi_{\beta(t)}}(t,0)\geq -
    \widetilde \chi_{\sR},
    $$
where
    \begin{equation}\label{chitildeRdef}
    \widetilde\chi_{\sR}=\inf
    \bigl\{ -\lambda_{\sR}(\psi)\colon \psi\in \Ccal^2(Q_{\sR}) \mbox{ and }
    {\mathcal L}_{\sR}(\psi)< 1 \bigr\}.
    \end{equation}
It remains to show that, for any $\rho>0$, we have
$\limsup_{R\uparrow\infty} \widetilde \chi_{\sR} \leq\widetilde \chi(\rho).$
This can be seen as follows: By Proposition~\ref{asconstant} the variational problem in
\eqref{chitildedef} has a minimizer $\psi^*$, a parabola with ${\mathcal L}(\psi^*)=1$.
Pick $\psi_{\sR}=\eps_{\sR}+ \psi^*|_{Q_{\sR}}$, where $\eps_{\sR}>0$ is chosen such that
${\mathcal L}_{\sR}(\psi_{\sR})=1-\sfrac 1R$. Obviously $\eps_{\sR}\downarrow 0$. It is easy to
show, using the explicit principal eigenfunction of $\Delta+\psi^*$ that
$\lim_{R\to\infty} \lambda_{\sR}(\psi_{\sR})=\lambda(\psi^*)$.
This completes the proof of the lower bound in (\ref{mainaimquen}) subject
to Proposition~\ref{microopt}.
\medskip

\noindent
We finally prove Proposition~\ref{microopt}:

\begin{Proof}{Proof of Proposition~\ref{microopt}.} This is very similar to the proof
of \cite[Prop.~5.1]{BK01}. Recall that $\psi_t(\cdot)=\alpha(\beta(t))^{-2}
\psi(\cdot\,\alpha(\beta(t))^{-1})$. Consider the event
$$    A^{\ssup{t}}_y=\bigcap_{z\in B_{R\alpha(\beta(t))}}\Big\{\xi_{\beta(t)}(y+z)\geq \psi_t(z)-
    \frac \eps{2\alpha(\beta(t))^2}\Big\},\qquad \mbox{ for } y\in \Z^d.$$
Note that the distribution of $A^{\ssup{t}}_y$ does not depend on $y$.
Our first goal is to show that, for every $\eps>0$,
    \begin{equation}\label{probAtesti}
    \Prob\big(A^{\ssup{t}}_0\big)\geq t^{-d\Lcal_{\sR}(\psi)-C\eps +o(1)},\qquad \mbox{ as }  t\uparrow\infty,
    \end{equation}
where $C>0$ depends only on $R$ and $\psi$, but not on $\eps$.

It is convenient to abbreviate
    \begin{equation}
    \label{st-def}
    s_t=\beta(t)\alpha(\beta(t))^{-d}.
    \end{equation}
Let $f\in\Ccal(Q_{\sR})$ be some positive auxiliary function (to be determined later), and consider the
tilted probability measure
    $$
    \Prob_{t,z}(\,\cdot\,)=\big\langle e^{f_t(z)\xi_{\beta(t)}(z)}\1\{\xi(z)\in\,\cdot\,\}
    \big\rangle e^{-H(f_t(z))+f_t(z)H(s_t)/s_t},\qquad \mbox{ for } z\in\Z^d,
    $$
where $f_t(z)=s_t f(z\alpha(\beta(t))^{-1})$ is the scaled version of $f$.
The purpose of this tilting is to make the event $A_0^{\ssup t}$ typical.
We denote the expectation with respect to $\Prob_{t,z}$ by $\langle\,\cdot\,\rangle_{t,z}$.
Consider the event
    $$
    D_t(z)=\Big\{\frac\eps{2\alpha(\beta(t))^2}\geq \xi_{\beta(t)}(z)-\psi_t(z)\geq-
    \frac\eps{2\alpha(\beta(t))^2}\Big\}.
    $$
Using that $\bigcap_{z\in B_{R\alpha(\beta(t))}}D_t(z)\subseteq A_0^{\ssup{t}}$ and the left inequality
in the definition of $D_t(z)$, we obtain
    \begin{equation}\label{microbox1}
    \begin{aligned}
    \Prob\big(A_0^{\ssup{t}}\big)&\geq \exp\Bigl\{\sum_{z\in B_{R\alpha(\beta(t))}}
    \Big[H(f_t(z))-f_t(z)\Big(\frac {H(s_t)}{s_t}+\psi_t(z)+\frac\eps{2\alpha(\beta(t))^2}\Big)\Big]\Big\}\\
    &\qquad\times\prod_{z\in B_{R\alpha(\beta(t))}}\Prob_{t,z}\big(D_t(z)\big).
    \end{aligned}
    \end{equation}
Since $\beta(t)\alpha(\beta(t))^{-2}=d\log t$, it is clear from a Riemann sum approximation  that
\begin{equation}\label{microbox2}
\begin{aligned}
\exp\Bigl\{\sum_{z\in B_{R\alpha(\beta(t))}} & \Big[-f_t(z)\Big(\psi_t(z)
+\frac\eps{2\alpha(\beta(t))^2}\Big)\Big]\Big\} \\
& = \exp\Big\{ - \frac{\beta(t)}{\alpha(\beta(t))^2}\, \frac{1}{\alpha(\beta(t))^d}\, 
\sum_{z\in B_{R\alpha(\beta(t))}} f\big( \sfrac{z}{\alpha(\beta(t)))} \big)
\Big( \psi\big( \sfrac{z}{\alpha(\beta(t)))} \big) + \sfrac\eps 2 \Big) \Big\} \\
& =t^{-d\langle f,\psi\rangle -d\frac \eps2\langle f,\1\rangle +o(1)},\qquad \mbox{ as } t\uparrow\infty.
\end{aligned}
\end{equation}    
We use the uniformity of the convergence in \eqref{basic}, the definitions
\eqref{alphadef} of $\alpha(\,\cdot\,)$  and \eqref{betatdef} of $\beta(t)$,
and a Riemann sum approximation to obtain
\begin{equation}\label{microbox3}
\begin{aligned}
\sum_{z\in B_{R\alpha(\beta(t))}}&\Big[H(f_t(z))-f_t(z)\sfrac {H(s_t)}{s_t}\Big]
=\sum_{z\in B_{R\alpha(\beta(t))}}\Big[H\Big(s_t f\big(\sfrac z{\alpha(\beta(t))}\big)\Big)
- f\big(\sfrac z{\alpha(\beta(t))}\big)H(s_t)\Big]\\
&=\kappa(s_t)\,\Big[ \sum_{z\in B_{R\alpha(\beta(t))}}\rho f\big(\sfrac z{\alpha(\beta(t))}\big)
\log f\big(\sfrac z{\alpha(\beta(t))}\big) + o(\alpha\circ\beta(t)^d) \Big] \\
&=\big(\Hcal_R(f)+o(1)\big)\,(1+o(1))\, \frac{\beta(t)}{\alpha(\beta(t))^2}=
\Hcal_R(f)\, \big(d (\log t) +o(1) \big).
\end{aligned}
\end{equation}
Using \eqref{microbox2} and \eqref{microbox3} in \eqref{microbox1}, we arrive at
$$    \Prob\big(A_0^{\ssup{t}}\big)\geq t^{d\big( \Hcal_{\sR}(f)-\langle f,\psi\rangle -
    \frac \eps2\langle f,\1\rangle \big) +o(1)}\prod_{z\in B_{R\alpha(\beta(t))}}\Prob_{t,z}\big(D_t(z)\big),
    \qquad \mbox{ as }t\uparrow\infty.$$
Recall from \eqref{calLRdef} that $\Lcal_{\sR}$ is the Legendre transform of $\Hcal_{\sR}$.
We choose $f$ as the minimizer on the right of  \eqref{calLRdef}, i.e., such that
$\Hcal_{\sR}(f)-\langle f,\psi\rangle=-\Lcal_{\sR}(\psi)$. Hence, to show that \eqref{probAtesti} holds,
it is sufficient to show that
    \begin{equation}\label{Desti1}
    \prod_{z\in B_{R\alpha(\beta(t))}}
    \Prob_{t,z}\big(D_t(z)\big)\geq t^{o(1)},\qquad \mbox{ as }t\uparrow\infty.
    \end{equation}
To show this, note that
    \begin{equation}\label{Desti2}
    \begin{aligned}
    \Prob_{t,z}\big(D_t(z)\big)&= 1-\Prob_{t,z}\Big\{\xi_{\beta(t)}(z)> \psi_t(z)+
    \frac\eps{2\alpha(\beta(t))^2}\Big\}\\
    &\qquad -\Prob_{t,z}\Big\{\xi_{\beta(t)}(z)< \psi_t(z)-\frac\eps{2\alpha(\beta(t))^2}\Big\}.
    \end{aligned}
    \end{equation}
Since both terms are handled in the same way, we treat only the second term. For any
$a>0$ we use the exponential Chebyshev inequality to bound
    $$
    \begin{aligned}
    \Prob_{t,z}&\Big\{\xi_{\beta(t)}(z)< \psi_t(z)-\frac\eps{2\alpha(\beta(t))^2}\Big\}\\
    &\leq e^{-H(f_t(z))+f_t(z)H(s_t)/s_t}\Big\langle
    \exp\Big\{f_t(z)\xi_{\beta(t)}(z)+a\Big[\psi_t(z)-\xi_{\beta(t)}(z)-
    \frac\eps{2\alpha(\beta(t))^2}\Big]\Big\}\Big\rangle\\
    &=e^{H(f_t(z)-a)-H(f_t(z))+aH(s_t)/s_t}e^{a[\psi_t(z)-\eps/2\alpha(\beta(t))^2]}.
    \end{aligned}
    $$
We pick $a=\delta_t f_t(z)$ with some $\delta_t\downarrow 0$.
Then the terms involving $H$ can be treated similarly to \eqref{microbox3}. Indeed,
abbreviating $\widetilde f=f(z\alpha(\beta(t))^{-1})$, we obtain
    $$
    \begin{aligned}
    H(&f_t(z)-a)-H(f_t(z))+a\frac{H(s_t)}{s_t}
    =H\big(s_t(1-\delta_t)\widetilde f \big)-(1-\delta_t)\widetilde f\, H\big(s_t\big)-
    \Big[H\big(s_t\widetilde f\big)-\widetilde f\, H\big(s_t\big)  \Big]\\
    &=\kappa(s_t)(1+o(1))\Big[\widehat H\big((1-\delta_t)\widetilde f\big)-\widehat
    H\big(\widetilde f\big)\Big]
    =-\frac{\rho+o(1)}{\alpha(\beta(t))^d}\widetilde f\,\delta_t\, d\,\log t
    \, \big[1+\log\widetilde f \big],
    \end{aligned}
    $$
where we also used the approximation $\log(1-\delta_t)=-\delta_t(1+o(1))$. Hence, we obtain
    $$
    \Prob_{t,z}\Big\{\xi_{\beta(t)}(z)< \psi_t(z)-\frac\eps{2\alpha(\beta(t))^2}\Big\}\leq
    t^{\delta_t\alpha(\beta(t))^{-d}\widetilde f\,[\widetilde \psi -\eps/2-\rho(1+\log\widetilde f)](d+o(1))},
    \qquad\mbox{ as } t\uparrow\infty,
    $$
where we recall \eqref{xibetadef} and abbreviate $\widetilde \psi=\psi(z\alpha(\beta(t))^{-1})$.
Recall that we chose $f$ optimally in \eqref{calLRdef}, which in particular means that $\log f(x)=\psi(x)/\rho-1$. Hence, for some $C>0$, not depending on $t$ nor on $z$, we have, for $t>1$ large enough,
    $$
    \Prob_{t,z}\Big\{\xi_{\beta(t)}(z)< \psi_t(z)-\frac\eps{2\alpha(\beta(t))^2}\Big\}
    \leq t^{-Cd\eps \delta_t\alpha(\beta(t))^{-d}}\leq \frac 14.
    $$
Going back to \eqref{Desti2} and assuming that the first probability term satisfies the same bound, we have
    \begin{equation}\label{microopt5}
    \begin{aligned}
    \prod_{z\in B_{R\alpha(\beta(t))}}\Prob_{t,z}\big(D_t(z)\big)
    \geq \big(1-\sfrac{1}{2}\big)^{|B_{R\alpha(\beta(t))}|}
    =e^{C(R\alpha(\beta(t)))^d}=e^{o(\log{t})}=t^{o(1)},
    \end{aligned}
    \end{equation}
where we use that $\alpha$ is slowly varying and $\beta(t)=(\log{t})^{1-o(1)}$,
so that $\alpha(\beta(t))^d \leq \beta(t)^{d\eta}=o(\log{t})$ for $t\to \infty$.
This proves \eqref{Desti1}, and therefore \eqref{probAtesti}.
\vskip0.2cm

We finally complete the proof of Proposition~\ref{microopt}. As in the proof of \cite[Prop.~5.1]{BK01}
it suffices to prove the almost sure existence
of a (random) $n_0\in\N$ such that, for any $n\geq n_0$, there is a $y_n\in B_{\gamma_{e^n}}$ such that
the event $A^{\ssup{e^{n+1}}}_{y_n}$ occurs. In the following, we abbreviate $t=e^n$. Let $M_t=B_{\gamma_t}
\cap \lfloor 3R\alpha(\beta(et))\rfloor \Z^d$. Note that $|M_t|\geq t^{d-o(1)}$ as $t\uparrow\infty$ and
that the events $A^{\ssup{et}}_y$ with $y\in M_t$ are independent. It suffices to show the summability of
    $$
    p_t=\Prob\Big\{\sum_{y\in M_t}\1\{A^{\ssup{et}}_y\}\leq \sfrac 12 |M_t|\Prob(A^{\ssup{et}}_0)\Big\}
    $$
on $t\in e^{\N}$. Indeed, since, by \eqref{probAtesti},
    \begin{equation}\label{summable}
    |M_t|\Prob(A^{\ssup{et}}_0)\geq t^{d-d\Lcal_{\sR}(\psi)-C\eps-o(1)}
    \end{equation}
tends to infinity if $\eps>0$ is small enough (recall that $\Lcal_{\sR}(\psi)<1$),
the summability ensures, via the Borel-Cantelli lemma, that, for all sufficiently
large $t$, even a growing number of the events $A^{\ssup{et}}_y$ with $y\in M_t$
occurs. To show the summability of $p_t$ for $t\in e^{\N}$, we use the
Chebyshev inequality to estimate
    $$
    p_t\leq \Prob\Big\{\Big[\sum_{y\in M_t}\1\{A^{\ssup{et}}_y\}-
    \Big\langle \sum_{y\in M_t}\1\{A^{\ssup{et}}_y\}\Big\rangle\Big]^2
    >\frac 14 \big[|M_t|\Prob(A^{\ssup{et}}_0)\big]^2\Big\} \leq
    4 \frac{1-\Prob(A^{\ssup{et}}_0)}{|M_t|\Prob(A^{\ssup{et}}_0)}.
    $$
The summability over all $t\in e^{\N}$ is clear from \eqref{summable}.
\qed\end{Proof}

\section{Appendix: Corrected proof of Lemma 4.2 in \cite{BK01}}\label{sec-erratum}

We use the opportunity to correct an error in the proof of one of the
main results of \cite{BK01}, the analogue of Theorem~\ref{momasy} for case~(4)
in Section~\ref{sec-uni}. In the original proof the large deviation principle of 
Proposition~\ref{prop-GKS04} and Varadhan's lemma are applied to the functional 
$f \mapsto - \int f^\gamma \, dx$, which fails to be continuous in the topology 
of the large deviation principle.
Here we adapt the techniques of the present paper to derive this result.
We use the notation of Section~\ref{sec-proofann}.

Recall case (4) from Section~\ref{sec-uni}. That is, we are in the case where $\esssup\,\xi(0)=0$,
$\gamma\in(0,1)$ and $\kappa^*=0$. The case $\gamma=0$ is easier and can be treated analogously.
The main assumption is that $\lim_{t\to\infty}
\widetilde H_t(x)=-D x^\gamma$,
uniformly in $x$ on compact subsets of $[0,\infty)$, where
    \begin{equation}
    \label{widetildeHdef}
    \widetilde H_t(x)=\frac{\alpha(t)^{d+2}}t H\Big(x\frac t{\alpha(t)^d}\Big),
    \end{equation}
and $D>0$ is a parameter. We have $\alpha(t)=t^{\nu+o(1)}$ as $t\to\infty$, where
$\nu=\frac{1-\gamma}{d+2-d\gamma}\in(0,\frac 1{d+2})$.

The step which needs amendment in \cite{BK01} is the following analogue of
Proposition~\ref{largedeviat1}:

\begin{prop}\label{prop-BKrepair}
~
\vskip-1cm

\begin{enumerate}
\item[(i)] For any $R>0$ and $M>0$, $\displaystyle{\limsup_{t\uparrow \infty}\frac {\alpha(t)^2}t\log \big\langle u^\xi_{R\alpha(t)}(t,0)\big\rangle\leq -\chi^{\ssup{M}}}$.

\item[(ii)] For any $R>0$, $\displaystyle{\liminf_{t\uparrow \infty}
\frac {\alpha(t)^2}t\log \big\langle u^\xi_{R\alpha(t)}(t,0)\big\rangle\geq -\chi_{\sR}}$,
\end{enumerate}
where
\begin{eqnarray}
\chi^{\ssup{M}}&=&\inf_{\heap{g\in H^1(\R^d)}{\|g\|_2=1}}\Big(\|\nabla g\|_2^2
+D\int \big(g^2(x)\land M\big)^\gamma\, dx\Big), \nonumber\\
\chi_{\sR}&=&\inf_{\heap{g\in H^1(\R^d)}{\|g\|_2=1,\supp(g)\subseteq Q_{\sR}}}\Big(\|\nabla g\|_2^2
+D\int g^{2\gamma}(x)\, dx\Big).\nonumber
\end{eqnarray}
\end{prop}

\begin{Proof}{Proof.} Introduce
    $$
    \Hcal^{\ssup t}_{\sR}(f)=\int_{Q_{\sR}}\widetilde H_t\big(f(x)\big)\, dx,\qquad \mbox{ for } f\in L^1(Q_{\sR}),
    \, f\geq 0.
    $$
As in \eqref{utident3}, we have
    \begin{equation}\label{errurepr}
    \big\langle u^\xi_{R\alpha(t)}(t,0)\big\rangle=
    \E_0\Big[\exp\Big(\frac t{\alpha(t)^2}\Hcal_{\sR}^{\ssup t}(L_t)\Big)\1\{\supp(L_t)\subseteq Q_{\sR}\}\Big],
    \end{equation}
where we recall the rescaled and normalized local times $L_t$.

We start with the proof of (i). Fix $M>0$. With $\Hcal_{\sR}(f)=-D\int_{Q_{\sR}} f(x)^\gamma\, dx$,
we have, uniformly in $f\in L^1(Q_{\sR})$, $f\geq 0$,
    $$
    \limsup_{t\uparrow\infty}\Hcal_{\sR}^{\ssup t}(f)\leq \limsup_{t\uparrow\infty}\Hcal_{\sR}^{\ssup t}
    (f\land M)=\Hcal_{\sR}(f\land M).
    $$
Note that $\Hcal_{\sR}(L_t\land M)=\alpha(t)^2 G_t(\sfrac 1t \ell_t)$, where we introduce
    $$
    G_t(\mu)=-\frac D{\alpha(t)^2}\alpha(t)^{-d}\sum_{z\in B_{R\alpha(t)}}\Big((\alpha(t)^d \mu(z))
    \land M\Big)^\gamma,\qquad \mbox{ for }\mu\in \Mcal_1(B_{R\alpha(t)}).
    $$
We now use Proposition \ref{prop-BHK04} for $B=B_{R\alpha(t)}$ and $F=G_t$ to obtain from
\eqref{errurepr} that, for any large $t$,
    $$
    \begin{aligned}
    \big\langle u^\xi_{R\alpha(t)}(t,0)\big\rangle&\leq e^{o(t\alpha(t)^{-2})}
    \E_0\Big[\exp\Big(t \, G_t(\sfrac 1t\ell_t)\Big)\1\{\supp(\ell_t)\subseteq B_{R\alpha(t)}\}\Big]\\
    &\leq e^{o(t\alpha(t)^{-2})}\exp\Big(-t  \chi^{\ssup{M}}_t\Big),
    \end{aligned}
    $$
where
    $$
    \chi^{\ssup{M}}_t=\inf_{\mu\in\Mcal_1(B_{R\alpha(t)})}\Big( \frac 12
    \sum_{x\sim y}\Big(\sqrt{\mu(x)}-\sqrt{\mu(y)}\Big)^2-G_t(\mu)\Big).
    $$
The proof of the upper bound is finished as soon as we have shown that
    \begin{equation}\label{erruppend}
    \liminf_{t\uparrow \infty}\alpha(t)^2\chi^{\ssup{M}}_t\geq \chi^{\ssup{M}}.
    \end{equation}
This is shown as follows. Let $(t_n \colon n\in\N)$ be a sequence of positive numbers $t_n\to\infty$
along which $\liminf_{t\uparrow \infty}\alpha(t)^2\chi^{\ssup{M}}_t$ is realized. We may
assume that its value is finite. Let $(\mu_n \colon n\in\N)$ be a sequence of approximative minimizers,
i.e., probability measures on $\Z^d$ having support in $B_{R\alpha(t)}$ such that
    $$
    \liminf_{n\to\infty}\Big[\alpha(t_n)^2\sfrac 1{2}\sum_{z\sim y}\Big(\sqrt{\mu_n(z)}-
    \sqrt{\mu_n(y)}\Big)^2+D\alpha(t_n)^{-d}\sum_z
    \Big(\big(\alpha(t_n)^{d}\mu_n(z)\big)\land M\Big)^\gamma\Big]
    $$
is equal to the left-hand side of \eqref{erruppend}. For any
$i\in\{1,\dots,d\}$ consider $g^{\ssup{i}}_n\colon \R^d\to\R$ given by
    $$
    \begin{aligned}
    g^{\ssup{i}}_n(x)&=\sqrt{\alpha(t_n)^d\mu_n\big(\lfloor \alpha(t_n)x\rfloor\big)}\\
    &\quad+\Big(x_i-\frac{\lfloor \alpha(t_n)x_i\rfloor}{\alpha(t_n)}\Big)\alpha(t_n)
    \Big(\sqrt{\alpha(t_n)^d\mu_n\big(\lfloor \alpha(t_n)x\rfloor+{\rm e}_i\big)}
    -\sqrt{\alpha(t_n)^d\mu_n\big(\lfloor \alpha(t_n)x\rfloor\big)}\Big),
    \end{aligned}
    $$
where ${\rm e}_i\in\Z^d$ is the $i^{\rm th}$ unit vector. For $x=(x_j \colon j=1,\dots,d)
\in\R^d$, we abbreviate $\widetilde x_i=(x_j \colon j\not=i)\in\R^{d-1}$ and denote
$g^{\ssup{i}}_{n,\widetilde x_i}(x_i)=g^{\ssup{i}}_n(x)$. For almost every
$\widetilde x_i\in\R^{d-1}$, the map $g^{\ssup{i}}_{n,\widetilde x_i}$ is
continuous and piecewise affine, and hence lies in $H^1(\R)$ with support in $[-R,R]$.
Furthermore,
    $$
    (g^{\ssup{i}}_{n,\widetilde x_i})'(x_i)=\frac{\partial g^{\ssup{i}}_n}{\partial x_i}(x)=
    \alpha(t_n)\Big(\sqrt{\alpha(t_n)^d
    \mu_n\big(\lfloor \alpha(t_n)x\rfloor+{\rm e}_i\big)}-\sqrt{\alpha(t_n)^d\mu_n
    \big(\lfloor \alpha(t_n)x\rfloor\big)}\Big).
    $$
Hence, using Fubini's theorem and Fatou's lemma, we see that
    $$
    \begin{aligned}
    \infty&>\liminf_{n\to\infty}\alpha(t_n)^2\frac 1{2}\sum_{z\sim y}\Big(
    \sqrt{\alpha(t_n)^d\mu_n(z)}-\sqrt{\alpha(t_n)^d\mu_n(y)}\Big)^2\\
    &=\liminf_{n\to\infty}\sum_{i=1}^d\int_{\R^{d-1}}d \widetilde x_i\int_{\R}d
    x_i \,\big|(g^{\ssup{i}}_{n,\widetilde x_i})'(x_i)\big|^2
    \geq \sum_{i=1}^d\int_{\R^{d-1}}d \widetilde x_i\,\liminf_{n\to\infty}\int_{\R}
    d x_i \,\big|(g^{\ssup{i}}_{n,\widetilde x_i})'(x_i)\big|^2.
    \end{aligned}
    $$
Since
    \begin{equation}
    \label{diffest}
    |x_i-\lfloor \alpha(t_n)x_i\rfloor/\alpha(t_n)|\leq \alpha(t_n)^{-1},
    \end{equation}
this also shows that
    \begin{equation}
    \label{convL2norm}
    \lim_{n\to\infty}\|g_n^{\ssup{i}}-\sqrt{\alpha(t_n)^d\mu_n(\lfloor
    \alpha(t_n)\,\cdot\,\rfloor)}\|_2=0.
    \end{equation}
In particular, $g_n^{\ssup{i}}$ is asymptotically
$L^2$-normalized. Furthermore, it follows that, along a suitable subsequence, for
almost all $\widetilde x_i\in\R^{d-1}$, $g^{\ssup{i}}_{n,\widetilde x_i}$ converges
to some $g^{\ssup{i}}_{\widetilde x_i}\in H^1(\R)$. The convergence is
(i)~strong in $L^2$, (ii)~pointwise almost everywhere, and (iii)~weak in $L^2$
for the gradients. The limit satisfies
    \begin{equation}
    \label{H1norm}
    \liminf_{n\to\infty}\alpha(t_n)^2\frac 1{2}\sum_{z\sim y}\Big(\sqrt{\alpha(t_n)^d
    \mu_n(z)}-\sqrt{\alpha(t_n)^d\mu_n(y)}\Big)^2\geq \sum_{i=1}^d\int_{\R^{d-1}}d
    \widetilde x_i\,\int_{\R}d x_i\, \big|(g^{\ssup{i}}_{\widetilde x_i})'(x_i)\big|^2.
    \end{equation}
Since $g^{\ssup{i}}_{n,\widetilde x_i}(x_i)=g^{\ssup{i}}_n(x)$ and
$\lim_{n\to\infty}\|g_n^{\ssup{i}}-\sqrt{\alpha(t_n)^d\mu_n(\lfloor \alpha(t_n)
\,\cdot\,\rfloor)}\|_2=0$, there is $g\in L^2(\R^d)$ such that $g(x)=
g^{\ssup{i}}_{\widetilde x_i}(x_i)$ for almost all $x\in\R^d$. In particular,
(a)~$g\in H^1(\R^d)$
with (b)~$\|g\|_2=1$, (c)~$\supp(g)\subset Q_{\sR}$ and (d)
    $$
    \|\nabla g\|_2^2\leq \liminf_{n\to\infty}\alpha(t_n)^2\frac 1{2}\sum_{z\sim y}
    \Big(\sqrt{\alpha(t_n)^d\mu_n(z)}-\sqrt{\alpha(t_n)^d\mu_n(y)}\Big)^2.
    $$
Indeed, (a) follows from~(b) and~(d). Item~(b) follows
from \eqref{convL2norm}, while item~(c) is trivially satisfied.
We are left to prove item~(d). Since $g^{\ssup{i}}_{\widetilde x_i}(x_i)=g(x)$
for almost every $x$, we get
    \eq
    (g^{\ssup{i}}_{\widetilde x_i})'(x_i)
    =\frac{\partial}{\partial x_i} g^{\ssup{i}}_{\widetilde x_i}(x_i)
    = \frac{\partial}{\partial x_i} g(x),
    \en
and hence
    \eq
    \sum_{i=1}^d\int_{\R^{d-1}}d
    \widetilde x_i\,\int_{\R}d x_i\, \big|(g^{\ssup{i}}_{\widetilde x_i})'(x_i)\big|^2
    =\int_{\R^{d}}d x \sum_{i=1}^d\big|\frac{\partial}{\partial x_i} g(x)\big|^2
    =\|\nabla g\|_2^2.
    \en
Therefore, item (d) follows from \eqref{H1norm}.

It remains to show that $\int (g(x)^{2}\land M)^\gamma\, dx\leq \liminf_{n\to\infty}
\alpha(t_n)^{-d}\sum_z ((\alpha(t_n)^d\mu_n(z))\land M)^\gamma$. Note that
    \begin{equation}\label{BKproof1}
    \begin{aligned}
    \alpha(t_n)^{-d}\sum_z & \Big(\big(\alpha(t_n)^d\mu_n(z)\big)\land M\Big)^\gamma
    =\int\Big(\big(\alpha(t_n)^d\mu_n(\lfloor \alpha(t_n)x\rfloor)\big)\land M\Big)^\gamma\, dx\\
    &=\int \Big(\Big|g_n^{\ssup{i}}(x)-\Big(x_i-\frac{\lfloor \alpha(t_n)x_i\rfloor}
    {\alpha(t_n)}\Big)(g^{\ssup{i}}_{n,\widetilde x_i})'(x_i)\Big|^{2\gamma}\land M^\gamma\Big)\, dx.
    \end{aligned}
    \end{equation}
We next use the inequality $|a-b|^{2\gamma}\geq (|a|^\gamma-|b|^\gamma)^2\geq
|a|^{2\gamma}-2|ab|^\gamma$ and for the subtracted term use Jensen's inequality
and the Cauchy-Schwarz inequality, as well as \eqref{diffest}, to see that
$$\begin{aligned}
    \int \Big|g_n^{\ssup{i}}(x) & \Big(x_i-\frac{\lfloor \alpha(t_n)x_i\rfloor}
    {\alpha(t_n)}\Big)(g^{\ssup{i}}_{n,\widetilde x_i})'(x_i)\Big|^{\gamma}\, dx\\
& \leq (2R)^{d} \Big[  (2R)^{-d} \int_{Q_R} \big|g_n^{\ssup{i}}(x) \big|
\Big| \Big(x_i-\frac{\lfloor \alpha(t_n)x_i\rfloor}
    {\alpha(t_n)}\Big)(g^{\ssup{i}}_{n,\widetilde x_i})'(x_i)\Big| \, dx \Big]^{\gamma} \\
& \leq \alpha(t_n)^{-\gamma} \, (2R)^{(1-\gamma)d} \;\big\|g_n^{\ssup{i}}\big\|_2^\gamma \;\;
\big\| \sfrac{\partial}{\partial x_i} g_n^{\ssup{i}} \big\|_2^\gamma,
\end{aligned}$$ 
which is negligible. Next we use the fact that $g_n^{\ssup{i}}\to g$ pointwise and
Fatou's lemma to see that the limit inferior of the right hand
side of \eqref{BKproof1} is not smaller than $\int (g(x)^{2\gamma}\land M^\gamma)\, dx$.
This completes the proof of \eqref{erruppend} and therefore the proof of (i).

We next turn to the proof of (ii). First we show that, for any $f\in\Ccal(Q_{\sR})$ and any
family of $L^1(Q_{\sR})$-normalized functions $f_t\in L^1(Q_{\sR})$
satisfying $f_t\to f$ in the weak topology induced
 by test integrals against all continuous functions,
\begin{equation}\label{errHconv}
\liminf_{t\uparrow\infty}\Hcal_{\sR}^{\ssup t}(f_t)\geq \Hcal_{\sR}(f).
\end{equation}
We fix a large $M>0$ and estimate $\Hcal_{\sR}^{\ssup t}(f_t)\geq\Hcal_{\sR}^{\ssup t}(f_t\land M)
+\Hcal_{\sR}^{\ssup t}(f_t\1\{f_t>M\})$. We first handle the first term.
Introduce $\phi(x)=x^\gamma$ and let $g_y(x)=(1-\gamma)y^\gamma+\gamma  y^{\gamma-1}x$
denote the tangent of $\phi$ at $y\in(0,\infty)$. By concavity, we have $\phi\leq g_y$ on $(0,\infty)$
 for any $y>0$. This implies that, as $t\uparrow \infty$, for any $\eps>0$,
$$
\begin{aligned}
\Hcal_{\sR}^{\ssup t}(f_t\land M)&=o(1)-D\int_{Q_{\sR}}\phi\big(f_t(x)\land M\big)\,dx
\geq o(1) -D\int_{Q_{\sR}} g_{f(x)\lor \eps}\big(f_t(x)\big)\, dx\\
&\geq o(1)-D(1-\gamma)\int_{Q_{\sR}} \big(f(x)\lor \eps\big)^\gamma\,
dx-D\gamma\int_{Q_{\sR}} f_t(x)\big(f(x)\lor \eps\big)^{\gamma-1}\, dx\\
&= o(1)-D (1-\gamma)\int_{Q_{\sR}} (f\lor\eps)^\gamma-D\gamma\int_{Q_{\sR}} f\, (f \lor \eps)^{\gamma-1},
\end{aligned}
$$
where in the last step we used that $(f\lor\eps)^{\gamma-1}$ is continuous and $f_t\to f$. Letting
$\eps\downarrow 0$, we see that $\liminf_{t\uparrow\infty}\Hcal_{\sR}^{\ssup t}(f_t\land M)
\geq \Hcal_{\sR}(f)$ for any $M>0$.

It remains to show that $\liminf_{M\uparrow \infty}\liminf_{t\uparrow\infty}\Hcal_{\sR}^{\ssup t}
(f_t\1\{f_t>M\})\geq 0$. Fix $\delta>0$ such that $\gamma+\delta<1$.
Recall \eqref{widetildeHdef}. Since $H$ is regularly varying with
exponent $\gamma$, by~\cite[Proposition 1.3.6]{BGT87},
there is an $M>0$ such that, for any sufficiently large $t$,
    $$
    \widetilde H_t(x)\geq- x^{\gamma+\delta},\qquad \mbox{for any }x>M.
    $$
Hence,
    $$
    \Hcal_{\sR}^{\ssup t}
    (f_t\1\{f_t>M\})\geq -\int_{Q_{\sR}} f_t(x)^{\gamma+\delta}\1\{f_t(x)>M\}\, dx
    \geq -M^{\gamma+\delta-1}\int_{Q_{\sR}} f_t(x)\, dx= -M^{\gamma+\delta-1},
    $$
since $f_t$ is $L^1$-normalized. This completes the proof of \eqref{errHconv}.

We complete the proof of Proposition \ref{prop-BKrepair}(ii) by
using \eqref{errHconv} in  \eqref{errurepr} and use
the lower bound of Varadhan's lemma in \cite[Lemma~4.3.4]{DZ98}
to conclude that the assertion in (ii) holds.
\qed\end{Proof}

{\bf Acknowledgements:} We would like to thank Laurens de Haan
for helpful discussions on regularly varying functions, and the
organisers of the Workshop on \emph{Interacting stochastic systems} in
Cologne, 2003, where this work was initiated. The work of the
first author was supported in part by Netherlands Organisation for
Scientific Research (NWO). The second author
would like to thank the German Science Foundation for awarding
a Heisenberg grant (realized in 2003/04), and the third author
would like to acknowledge the support of the Nuffield Foundation
through grant NAL/00631/G and the EPSRC through grant EP/C500229/1.


\end{document}